\documentclass[11pt]{article} 
\usepackage{multicol} 
\usepackage{newlfont}
\usepackage{amsmath,amssymb}
\usepackage{tikz}
\usepackage{tikz-cd}

\def\rit{\mathbb{R}}
\def\zit{\mathbb{Z}}

\def\ppit{\mathbb{P}} 
\def\qit{\mathbb{Q}} 
\def\cit{\mathbb{C}}

\newcommand{\pf}{{\em Proof.~}}
\newcommand{\qed}{\hfill~~\mbox{$\Box$}}

\newenvironment{proof}{\smallskip \noindent \pf}{\qed \bigskip}

\newtheorem{theorem}{Theorem}[section]
\newtheorem{proposition}[theorem]{Proposition}
\newtheorem{definition}[theorem]{Definition}
\newtheorem{lemma}[theorem]{Lemma}
\newtheorem{corollary}[theorem]{Corollary}

\newtheorem{remark}[theorem]{Remark}
\newtheorem{example}[theorem]{Example}

\DeclareMathOperator{\card}{Card}

\DeclareMathOperator{\boite}{Box}

\DeclareMathOperator{\gr}{gr}
\DeclareMathOperator{\Spec}{Spec}

\DeclareMathOperator{\orb}{orb}
\DeclareMathOperator{\Card}{Card}
\DeclareMathOperator{\supp}{supp}
\DeclareMathOperator{\age}{age}

\DeclareMathOperator{\prim}{prim}
\DeclareMathOperator{\MHS}{MHS}
\DeclareMathOperator{\PF}{PF}

\begin{document}



\title{\bf Mixed Hodge structures for vanishing cycles and orbifold cohomology}
\author{\sc Antoine Douai \\
Universit\'e C\^ote d'Azur, CNRS, LJAD, FRANCE\\
Email address: antoine.douai@univ-cotedazur.fr}

\maketitle

\begin{abstract}
Above a Laurent polynomial $f$ one makes grow 
a vector space of vanishing cycles (after the work of Sabbah, singularity setting), a graded Milnor ring (after the work of Kouchnirenko) and an orbifold cohomology ring (after the work of Borisov, Chen and Smith). Under suitable assumptions, these structures are isomorphic and these identifications are interesting because some results are more explicit in one setting than in another.
 In particular, and in order to understand better the real structures and the dualities appearing in the singularity setting, we first look for the counterpart of Sabbah's mixed Hodge structures, initially defined on the space of vanishing cycles, on the orbifold cohomology ring. Then, we discuss to what extent the orbifold Poincar\'e duality defined by Chen and Ruan provides a polarization of this mixed Hodge structure.
We study in details the Hodge-Tate case, which can be read off from the ages of the sectors, a variation of the hard Lefschetz condition introduced by Fernandez. 
These notes go along with prior works of Fernandez and Wang.
\end{abstract}

\section{Introduction}
Let $f$ be a convenient and nondegenerate Laurent polynomial on $(\cit ^*)^n$ in the sense of Kouchnirenko \cite{K}. One associates with $f$:
\begin{itemize}
\item a mixed Hodge structure $\MHS_f =(H, F^{\bullet}, W_{\bullet})$ (that we call Sabbah's mixed Hodge structure, see \cite{Sab0}, \cite{Sab}) where $H$ denotes the space of vanishing cycles; the weight filtration $W_{\bullet}$ is the weight filtration of a nilpotent endomorphism of $H$ (the logarithm of the unipotent part of the monodromy at infinity of $f$) and the Hodge filtration $F^{\bullet}$ has a very simple description via the Brieskorn lattice of $f$ and a suitable Kashiwara-Malgrange $V$-filtration,
\item the $\qit$-graded Milnor ring 
$A^* := \gr^{\mathcal{N}}_{*} (\cit [u,u^{-1}]/(u\partial f /\partial u))$,
where $\mathcal{N}_{\bullet}$ denotes a Newton filtration suitably normalized and $(u\partial f/\partial u)$ denotes the ideal generated by the partial derivatives $(u_1 \partial f/\partial u_1,\ldots , u_n \partial f/\partial u_n)$ of $f$ (the basic reference for this framework is Kouchnirenko's paper \cite{K}),
\item if the Newton polytope $P$ of $f$ is simplicial,
the orbifold cohomology ring $H^{2*}_{\orb} (\mathcal{X}_P):=H^{2*}_{\orb} (\mathcal{X}_P, \cit )$ equipped with the orbifold cup-product, of a suitable toric Deligne-Mumford stack $\mathcal{X}_P$ in the sense of \cite{BCS} (by {\em loc. cit.}, a toric Deligne-Mumford stack corresponds to 
a stacky fan, that is a simplicial fan with a distinguished lattice point on each ray of the fan; 
in order to define $\mathcal{X}_P$, we consider the stacky fan whose simplicial fan $\Sigma_P$ is the one over the faces of $P$ and whose lattice points on the rays are the vertices of the polytope $P$).
\end{itemize}
Sabbah's mixed Hodge structures have some interesting applications, which justify the interest shown in it: see for instance \cite{D9}, \cite{DoMa}, \cite{DoSa1}, \cite{DoSa2}, \cite{Mann1} for the relations with Frobenius manifolds and mirror symmetry. By \cite[Theorem 4.3]{Sab0}, it is also isomorphic to the limit mixed Hodge structure constructed by Steenbrink and Zucker in \cite{SZ} on $H^{\bullet}(U, f^{-1} (t); \qit)$ when $t\rightarrow \infty$ (we set here $U:= (\cit^*)^n$).

The first part of this text surveys the relationships between these three settings.
 On the one hand, if $f(u)=\sum_{b\in \mathcal{V}(P)} u^b$
where $\mathcal{V} (P)$ denotes the set of the vertices of the simplicial Newton polytope $P$ (in what follows, we always assume that $f$ has this form),
we have 
an isomorphism of graded rings 
\begin{equation}\label{eq:MirIntro}
A^* \cong H^{2*}_{\orb} (\mathcal{X}_P)
\end{equation}
where $H^{2*}_{\orb} (\mathcal{X}_P)$ is equipped with the orbifold cup-product, 
see \cite[Theorem 1.1]{BCS}.
On the other hand, and by the very definition of the Hodge filtration of the mixed Hodge structure $\MHS_f$, we have an isomorphism
$$\gr_F^p H \cong \bigoplus_{n-1-p <\alpha \leq n-p} A^{\alpha}$$
for $p\in\zit$. 
 Using the characteristic properties of the mixed Hodge structures, we finally get 
isomorphisms
\begin{equation}\label{eq:MirrorTotal}
H\stackrel{\sim}{\longrightarrow} A^* \stackrel{\sim}{\longrightarrow} H^{2*}_{\orb} (\mathcal{X}_P ), 
\end{equation}
see Section \ref{sec:Resume}.
These isomorphisms are useful because some structures are described more easily in one setting than in another: 
the ring structure is obviously simpler on $A^*$, while the linear structure is more explicit on  
$H^{2*}_{\orb} (\mathcal{X}_P )$, thanks to the descriptions given in \cite{BCS} with the help of Stanley-Reisner rings. 
The main motivation of these notes was to use these isomorphisms in order to understand better the real structure and the dualities appearing in the singularity theory setting, and basically defined by the Lefschetz thimbles and K. Saito's higher residue pairings \cite{Ph}, \cite{SaiK}, by constructing explicit counterparts in the orbifold cohomology setting.

We first look for the image of the mixed Hodge structure $\MHS_f$ in $H^{2*}_{\orb} (\mathcal{X}_P )$.
Using the direct sum decomposition
\begin{equation}\label{eq:DecompOrbi}
H_{\orb}^{2\alpha } (\mathcal{X}_P) :=\bigoplus_{\ell \in F} H^{2(\alpha -i_{\ell})}(X_{\ell})
\end{equation}
where $i_{\ell}$ denotes the age of the sector $X_{\ell}$ and 
$F$ labels the set of the sectors, we define in Section \ref{sec:Filtrations} a Hodge filtration and a weight filtrations on the orbifold cohomology, which correspond to the filtrations of the mixed Hodge structure $\MHS_f$ via the isomorphisms (\ref{eq:MirrorTotal}). 
Originally, these filtrations where inspired by the description given in \cite{DoSa2} of the mixed Hodge structure $\MHS_f$ for mirror partners of weighted projective spaces and match with the ones defined in a combinatorial setting by Wang \cite{Wa}.

The last ingredient that we need in order to define a mixed Hodge structure is a conjugation on $H_{\orb}^{2*} (\mathcal{X})$. It turns out that, once the Hodge filtration and the weight filtration are defined, the correct ones to consider on $H_{\orb}^{2*} (\mathcal{X})$ are those which satisfy
\begin{equation}\nonumber 
\overline{H^{2i} (X_{\ell})}\subset \oplus_{j\geq i} H^{2j} (X_{\ell^{-1}})
\end{equation}
for each sector $X_{\ell}$,
where $\ell^{-1}$ denotes the image of $\ell$ under a natural involution on the labels (see Section \ref{sec:Conjug}).
We will call such conjugations {\em Poincar\'e-Lefschetz conjugations} (the conjugations of Sabbah's mixed Hodge structures provide examples of such Poincar\'e-Lefschetz conjugations, see for instance \cite[Section 5]{DoSa2} and Lemma \ref{lemma:ConjugWPSDoSa} in the appendix). We show in Theorem \ref{theo:MHSOrbifold}, using the correspondence between mixed Hodge structure and bigradings stated in \cite[Theorem 2.13]{CKS}, that each Poincar\'e-Lefschetz conjugation gives rise to a mixed Hodge structure, denoted by $\MHS_{\orb}$, on $H_{\orb}^{2*} (\mathcal{X})$.

On the way, we compute the Hodge numbers of the mixed Hodge structures $\MHS_{\orb}$ and $\MHS_f$
 in terms of the Betti numbers of the sectors $X_{\ell}$, see Proposition \ref{prop:GraduesOrbi}. Let us point out that, keeping in mind the identification between Sabbah's mixed Hodge structure and the limit mixed Hodge structure 
on $H^{\bullet}(U, f^{-1} (t); \qit)$ alluded to above, 
the computation of such Hodge numbers is also considered in other papers using different techniques (basically, motivic Milnor fibers instead of Brieskorn lattice), see for instance \cite{KS}, \cite{MT} and the survey \cite{Tak}. The calculation of (orbifold, irregular) Hodge numbers is a central issue in several aspects of mirror symmetry and we refer to \cite{GKR}, \cite{HL} and \cite{KKP} for other perspectives.

It follows immediately that the mixed Hodge structure $\MHS_{\orb}$ is of Hodge-Tate type if and only if 
\begin{equation}\label{eq:HTIntro}
[i_{\ell}]= [i_{\ell^{-1}}]\ \mbox{for all}\ \ell\in F
\end{equation}
(an analogous statement is shown in \cite{D13}, using a different method). Naturally enough, we will call (\ref{eq:HTIntro}) the {\em Hodge-Tate condition} (when the orbifold cohomology ring is integrally graded, this condition is already considered by Fernandez in \cite[(4.2)]{Fe}). 
We analyze some of its consequences in Section \ref{sec:HThpq}. Note that this condition is not really common and, after \cite[Proposition 3.3]{Sab2}, this produces numerous counter-examples to a conjecture by Kontsevitch, Kapranov and Pantev \cite{KKP} (see \cite{D13} for details).

In a different combinatorial framework, a similar study is carried out by Wang in \cite{Wa} using the mixed Hodge structure defined on the ring of conewise linear functions by McMullen \cite{McM} (see also \cite{FK}). Our filtrations are analogous to the ones considered in loc. cit.
The corresponding formulas for the Hodge numbers, although stated in different mathematical settings (Wang's formula is expressed in terms of the $h$-vector of the polytope $P$), complement each other: depending on the context, one formula could be applied rather than another. However, it should be emphasized that the mixed Hodge structures $\MHS_{\orb}$ are not defined a priori as a direct sum of pre-existing mixed Hodge structures on the cohomology of each sector.

If one looks for a polarization of the mixed Hodge structure $\MHS_{\orb}$, it is natural to use the orbifold Poincar\'e duality defined by Chen and Ruan in \cite{CR} since it shares many properties with the duality of $\mathcal{D}$-modules studied in \cite{Sab} (this is highlighted in Appendix \ref{sec:DualityDmodules} and Appendix \ref{sec:DualityDmodulesWPS}). 
 In Theorem \ref{theo:Polar} we get such a polarization by considering a suitable normalization/Tate-twist of this orbifold Poincar\'e duality. 
Of course, it is possible to define different polarizations, depending on the specific Poincar\'e-Lefschetz conjugations involved: see Section \ref{sec:Alternate} where another polarization is discussed, close to the one constructed by Fernandez in \cite[Theorem 5.3]{Fe} (in this last case we have to assume moreover that the Hodge-Tate condition (\ref{eq:HTIntro}) is satisfied).

Throughout this text, we will test our results on weighted projective spaces (the corresponding polytopes are reduced simplices in the sense of \cite{Conrads}), see for instance \cite{CCLT}, \cite{DoMa}, \cite{Mann}, \cite{Mann1} (mirror partners on the singularity side appear already in \cite{DoSa2}).
This is a nice setting to work with: quick computations can be carried out and, remarkably, almost everything can be deduced from the weight of the corresponding simplex. For instance, we get a very simple and effective formula for the Hodge numbers of the corresponding mixed Hodge structures $\MHS_{\orb}$ and $\MHS_f$, depending only on the ages of the sectors (see Corollary \ref{coro:NbreHodgeSimplex}). We get also nice formulas for the primitive subspaces and their Hodge decomposition 
and this allows us to get a polarization on the mixed Hodge structure $\MHS_{\orb}$ rather inexpensively, see Proposition \ref{prop:PolarWPS}.

These notes are organized as follows: in Section \ref{sec:Setting}, we describe and compare the three different settings considered in this paper; in Section \ref{sec:MixedHodgeOrbi}, we define a mixed Hodge structure on the orbifold cohomology and we give a formula for its Hodge numbers; we analyze the Hodge-Tate condition in Section \ref{sec:HThpq} and we discuss the existence of a polarization via the orbifold Poincar\'e duality in Section \ref{sec:Polarization}. The definitions about polarized mixed Hodge structures used in this text are recalled in Appendix \ref{sec:PolarizedMHS}. The link between the polarization constructed in this paper and the duality of $\mathcal{D}$-modules is emphasized in Appendix \ref{sec:DualityDmodules} and Appendix \ref{sec:DualityDmodulesWPS}.

\section{The settings}

\label{sec:Setting}

This section is expository in nature.

\subsection{The graded Milnor ring}
\label{sec:GradedMilnorRing}
Let $f(u) =\sum_m a_m u^m$, where $m=(m_1 ,\ldots , m_n)\in\zit^n$ and $u^m= u_1^{m_1}\ldots u_n^{m_n}$, be a Laurent polynomial defined on $(\cit^*)^n$. Let $\supp (f)=\{m , a_m \neq 0\}$. 
The {\em Newton polytope} $P$ of $f$ is the convex hull of $\supp (f)$ in $\rit^n$. We will always assume that $f$ is convenient (i.e. $P$ contains the origin as an interior point) and nondegenerate in the sense of Kouchnirenko \cite{K}.

Let $\mathcal{A}_f :=\mathcal{B}/ \mathcal{L}$ where 
$\mathcal{B}:=\cit [u_1 , u_1^{-1},\ldots , u_n ,u_n^{-1}]$
and $\mathcal{L}:=(u_1\partial f / \partial u_1,\ldots , u_n \partial f/ \partial u_n)$
is the ideal generated by the partial derivatives
 $u_1\partial f / \partial u_1 ,\ldots , u_n\partial f/ \partial u_n$ of $f$.
We define an increasing filtration $\mathcal{N}_{\bullet}$ on $\mathcal{B}$, indexed by $\qit$, by setting  
\begin{equation}\nonumber
\mathcal{N}_{\alpha}\mathcal{B}:=\{g\in \mathcal{B},\ \supp (g) \in \nu^{-1}(]-\infty ; \alpha ]) \}
\end{equation} 
where $\nu :\rit^n \rightarrow \rit$ is the Newton function of the Newton polytope $P$ of $f$ (the function which takes the value $1$ at the vertices of $P$ and which is linear on each cone generated by the faces of $P$) and
$\supp (g)=\{m\in\zit^n,\ a_m \neq 0\}$ if $g=\sum_{m\in \zit^n}a_m u^m \in\mathcal{B}$. 
By projection, the Newton filtration $\mathcal{N}_{\bullet}$ on $\mathcal{B}$ induces a filtration, also denoted by $\mathcal{N}_{\bullet}$, on $\mathcal{A}_f$. We will denote by $A^*$ the graded ring 
\begin{equation}\nonumber 
A^* := \gr^{\mathcal{N}}_{*} \mathcal{A}_f 
\end{equation}
and we call it the {\em graded Milnor ring}. Note that the multiplication by $f$ on 
$\mathcal{A}_f$ induces a map $[f]$ of degree one on $A^*$. 
This is Kouchnirenko's setting \cite{K}.

\subsection{Orbifold cohomology and simplicial polytopes}

\label{sec:BCSsetting}

We now briefly recall the setting of Borisov, Chen and Smith \cite{BCS}.
Let $P$ be a simplicial lattice polytope in $\rit^n$, containing the origin as an interior point. We will denote by $\mathcal{V}(P)$ the set of its vertices. 
Let
\begin{itemize}
\item $\mathcal{X}_{P}$ be the Deligne-Mumford stack associated by \cite[Section 3]{BCS} with the stacky fan $\mathbf{\Sigma}_{P} :=(\zit^n, \Sigma_{P} , \mathcal{V}(P))$, where $\Sigma_P$ is the fan whose cones are the cones over the faces of $P$,
\item $I_{\mathcal{X}_{P}}=\coprod_{\ell\in F} \mathcal{X_{\ell}}$ be the decomposition into connected components of the inertia orbifold of $\mathcal{X}_{P}$,
\item $H_{\orb}^{2\alpha}(\mathcal{X}_{P}; \cit ):=\oplus_{\ell\in F} H^{2(\alpha-i_{\ell})}(X_{\ell}; \cit )$ be the orbifold cohomology groups of $\mathcal{X}_{P}$ where $X_{\ell}$ denotes the underlying space of the sector $\mathcal{X}_{\ell}$ and $i_{\ell}:=\age (\mathcal{X}_{\ell})$ denotes its age. 
\end{itemize}
To simplify the notations, we will write $\mathcal{X}$ instead of $\mathcal{X}_P$, $H_{\orb}^{2\alpha}(\mathcal{X} )$ instead of $H_{\orb}^{2\alpha}(\mathcal{X}_{P}, \cit )$ and 
$H^{2k}(X_{\ell})$ instead of $H^{2k}(X_{\ell}; \cit )$. 

We will consider the decomposition (of $\cit$-vector spaces)
\begin{equation}\label{eq:NotationDecompHorb}
H_{\orb}^{2*} (\mathcal{X})=H_{\orb}^{2*} (\mathcal{X})_0 \oplus H_{\orb}^{2*} (\mathcal{X})_{\neq 0}
\end{equation}
where
$$H_{\orb}^{2\alpha } (\mathcal{X} )_0 :=\bigoplus_{\ell \in F,\ i_{\ell}\in\zit} H^{2(\alpha -i_{\ell})}(X_{\ell})\ \mbox{and}\ 
H_{\orb}^{2\alpha } (\mathcal{X})_{\neq 0} :=\bigoplus_{\ell \in F,\ i_{\ell}\notin\zit} H^{2(\alpha -i_{\ell})}(X_{\ell}).$$
We will also make a repeated use of the equality
\begin{equation}\label{eq:FormuleDim}
n_{\ell}+i_{\ell}+i_{\ell^{-1}}=n
\end{equation}
given by \cite[Lemma 3.2.1]{CR}, where $n_{\ell}:=\dim X_{\ell}$ and $\ell^{-1}$ denotes the image of $\ell$ under the involution of $F$ induced by the involution $I$ of $I_{\mathcal{X}}$ defined by $I( (p, (g))=(p, (g^{-1}))$ where $(g)$ denotes the conjugacy class of $g$ in $G_p$ (the local group at $p$), see \cite[Remark 3.1.4]{CR}.

 In our toric setting, we have the following concrete descriptions: 

\begin{itemize}
\item for $\sigma$ a cone in the fan $\Sigma_P$, let
$\boite (\sigma )$ be the set of the elements $v\in N$ such that 
$$v=\sum_{\rho_i \subseteq\sigma} q_i b_i \  \mbox{for some}\ 0< q_i <1$$
 where $\rho_i$ is the ray generated by the vertex $b_i$ of $P$. Let $\boite (\mathbf{\Sigma}_P )$ be the union of $\boite (\sigma )$ for all cones $\sigma \in\Sigma_P$. By \cite[Proposition 4.7]{BCS}, the sectors are parametrized by $v \in\boite (\mathbf{\Sigma}_P )$.
\item For $v \in\boite (\mathbf{\Sigma}_P )$, $X_{v}$ is the simplicial toric variety of the quotient fan $\Sigma_P /\sigma (v)$, where $\sigma (v)$ denotes the smallest cone of $\Sigma_P$ containing $v$ and $n_v:=\dim X_{v}=n-\dim \sigma (v)$. 
\item If $v\in\boite (\mathbf{\Sigma}_P )$, the age of the sector $X_v$ is $i_{v}=\nu (v)$ where $\nu $ is the Newton function of $P$ defined in Section \ref{sec:GradedMilnorRing}, see \cite[Remark 5.4]{BCS}.
\item For $v=\sum_{\rho_i\subseteq \sigma (v)} q_{i} b_i\in\boite (\sigma (v))$, let us define 
$v^{-1}:=\sum_{\rho_i\subseteq \sigma (v)}(1-q_{i}) b_{i}$.
Then $v^{-1}\in\boite (\sigma (v))$ and 
\begin{equation}\nonumber 
\nu (v) +\nu (v^{-1})=\dim \sigma (v)=n-n_v
\end{equation}
which is precisely equality (\ref{eq:FormuleDim}) above. The involution $I$ on the labels is defined by $I(v)=v^{-1}$. 
\end{itemize}

\subsection{Sabbah's mixed Hodge structure}

\label{sec:SabbahMHS}

In this section, we recall Sabbah's setting \cite{Sab0}, \cite{Sab}. 
Let $f$ be a convenient and nondegenerate Laurent polynomial on $U :=(\cit^*)^n$, where $n\geq 2$. We will denote by $G$ the localized Laplace transform of its Gauss-Manin system and by $G_0$ its Brieskorn module: by the very definition, we have
\begin{equation}\nonumber
G= \Omega^n (U) [\theta , \theta^{-1}]/(\theta d-df\wedge ) \Omega^{n-1} (U) [\theta , \theta^{-1}]
\end{equation}
and $G_0$ is the image of  $\Omega^n (U) [\theta ]$ in $G$.
Since $f$ has only isolated singularities (see \cite{K}), we have
\begin{equation}\nonumber
G_0 = \Omega^n (U) [\theta ]/(\theta d-df\wedge ) \Omega^{n-1} (U) [\theta ]
\end{equation}
and 
\begin{equation}\nonumber
G_0 /\theta G_0 \cong \Omega^{n}(U)/df\wedge \Omega^{n-1} (U),
\end{equation}
the latter being a finite dimensional vector space of dimension $\mu$, the global Milnor number of $f$. 
Actually, $G_0$ is a lattice in $G$, i.e $G_0$ is free over $\cit [\theta]$ and $\cit [\theta , \theta^{-1}]\otimes_{\cit [\theta ]} G_0 =G$: this is shown in \cite[Section 4]{DoSa1}, using the main properties of convenient and nondegenerate functions established in \cite{K}. This is a characteristic property of tame functions, for which Sabbah's mixed Hodge structures can be constructed \cite[Section 13]{Sab}.

Let $\tau :=\theta^{-1}$. Then $G$ is a free $\cit [\tau ,\tau^{-1}]$-module, equipped with a derivation $\partial_{\tau}$ (induced from the formula $\partial_{\tau} [\omega ]=[-f\omega]$ if $\omega\in\Omega^{n} (U)$ where the brackets denote the class in $G$).
Let $V_{\bullet}G$ be the Kashiwara-Malgrange filtration of $G$ along $\tau =0$ as in \cite[Section 1]{Sab}. This is an increasing filtration, indexed by $\qit$, which satisfies the following properties:
\begin{itemize}
\item for every $\alpha\in\qit$, $V_{\alpha} G$ is a free $\cit [\tau]$-module of rank $\mu$,
\item $\tau V_{\alpha} G \subset V_{\alpha -1} G$ and $\partial_{\tau} V_{\alpha} G \subset V_{\alpha +1} G$,
\item $\tau\partial_{\tau} +\alpha$ is nilpotent on $\gr_{\alpha}^V G$.
\end{itemize} 
We define, for $\alpha\in\qit$,
\begin{equation}\nonumber
H_{\alpha} := \gr_{\alpha}^V G\ \mbox{and}\ H:=\oplus_{\alpha\in [0,1[} H_{\alpha}
\end{equation}
and we write $H = H_{\neq 0}\oplus H_0$,
where $H_{\neq 0}:= \oplus_{\alpha \in ]0,1[} H_{\alpha}$.
The $\cit$-vector space $H$ (of dimension $\mu $) is equipped with a nilpotent endomorphism $N$ induced by $-(\tau\partial_{\tau}+\alpha )$ on $H_{\alpha}$.\\

\noindent {\bf  {\em The Hodge filtration.}} Let $G_{\bullet}$ be the increasing filtration of $G$ defined by $G_p :=\theta^{-p} G_0 =\tau^p G_0$ and
\begin{equation}\nonumber
G_p H_{\alpha} := G_p \cap V_{\alpha} G/G_p \cap V_{<\alpha} G
\end{equation}
for $\alpha\in [0,1[$ and $p\in\zit$. 
 The {\em Hodge filtration} on $H$ is the decreasing filtration $F^{\bullet}$ defined by
\begin{equation} \label{eq:HodgeFiltration}
F^p H_{\alpha} := G_{n-1-p}H_{\alpha}\ \mbox{if}\ \alpha\in ]0,1[\ \mbox{and}\ F^p H_{0} := G_{n-p}H_{0}
\end{equation}
for $p\in\zit$.\\

\noindent {\bf {\em The real structure and the weight filtration.}}
We have on $H$ a real structure coming from the identification 
\begin{equation}\nonumber
H\stackrel{\sim}{\longrightarrow} H^n (U, f^{-1}(t) ;\cit )= H^n (U, f^{-1}(t) ;\rit )\otimes_{\rit}\cit
\end{equation}
(for $|t|>>0$): in this setting, $H_{\alpha}$ corresponds to the generalized eigenspace of the monodromy at infinity associated with the eigenvalue $\exp (2i\pi\alpha )$ and the unipotent part of this monodromy is equal to $\exp (2i\pi N)$. We will write $H=H_{\rit}\otimes\cit$.
 We have the decomposition $H_{\rit}=H_{\rit ,0}\oplus H_{\rit, \neq 0}$ 
and the {\em weight filtration} $W_{\bullet}$ on $H_{\rit}$ is the increasing weight filtration of $2i\pi N$ centered at $n-1$ on $H_{\rit\neq 0}$ and centered at $n$ on $H_{\rit ,0}$.\\

\noindent {\bf {\em Sabbah's mixed Hodge structures.}} By \cite[Theorem 13.1]{Sab},
the tuple $\MHS_f :=(H, H_{\rit}, F^{\bullet}, W_{\bullet} )$ is a mixed Hodge structure.

\begin{remark} \label{rem:NGrF}
Let $\MHS_f$ be a mixed Hodge structure as above.
Deligne's subspaces $I^{p,q}$ (see \cite{Del} and Section \ref{sec:MHS} below) provide
the splittings $F^{p} H =\bigoplus_{i\geq p} I^{i,j}$ and $W_m =\bigoplus_{i+j\leq m} I^{i,j}$. Using the decomposition $F^{p} H =(\oplus_{j} I^{p,j})\bigoplus F^{p+1}H$, we get an isomorphism $H\cong\gr_F H$ such that $F^p H \cong  \bigoplus_{i\geq p} \gr_F^i H$ and, since $N(I^{p,q})\subset I^{p-1, q-1}$, a commutative diagram
$$\begin{tikzcd}
H \arrow[rightarrow] {r}{\sim} \arrow{d}{N} &  \gr_F H\arrow{d}{[N]}\\
H \arrow[rightarrow] {r}{\sim} & \gr_F H 
\end{tikzcd}$$
where $[N]$ is induced by $[N] :\gr_F^{i} H\rightarrow \gr_F^{i-1} H$: in particular, $N$ and $[N]$ have the same Jordan normal form (compare with \cite[Theorem 7.1]{ScSt}).
\end{remark}

\subsection{Isomorphisms} 
\label{sec:Iso}

We recall the mirror theorem of Borisov, Chen and Smith \cite{BCS}.

\begin{proposition}
[\cite{BCS}]
\label{prop:Mir} 
Let $P$ be a simplicial lattice polytope in $\rit^n$ containing the origin as an interior point and
let $f$ be the Laurent polynomial defined on $(\cit^*)^n$ by  
\begin{equation}\label{eq:fOrbi}
f (u) =\sum_{b\in \mathcal{V}(P)} u^b
\end{equation}
where $\mathcal{V}(P)$ denotes the set of the vertices of its Newton polytope $P$. Then, we have an isomorphism
\begin{equation}\label{eq:Mir}
\varphi : H_{\orb}^{2*}(\mathcal{X} )\stackrel{\sim}{\longrightarrow} A^*
\end{equation}
of $\qit$-graded rings where $H_{\orb}^{2*}(\mathcal{X} )$ is the orbifold cohomology of $P$, equipped with the orbifold cup-product, and $A^*$ is the graded Milnor ring of $f$.
\end{proposition}
\begin{proof}
Note first that the assumptions of Proposition \ref{prop:Mir} imply that $f$ is convenient and nondegenerate.
It remains to identify the graded ring appearing in \cite[Theorem 1.1]{BCS} with $A^*$. Recall that if $\partial : \mathcal{C}\rightarrow \mathcal{A}$ is a strict morphism of filtered module (i.e $\partial (\mathcal{C})\cap \mathcal{A}_q=
\partial (\mathcal{C}_q)$) we have $\gr (\mathcal{A}/ \partial \mathcal{C})\cong \gr (\mathcal{A})/ \gr (\partial) \gr (\mathcal{C})$, see for instance \cite[Lemme 4.4]{K}. 
We get the expected identification by applying this observation to the morphism
\begin{equation}\nonumber
\partial :\mathcal{A}^n \longrightarrow \mathcal{A}
\end{equation}
defined by 
\begin{equation}\nonumber
\partial (g_1 ,\ldots , g_n )= g_1 u_1 \frac{\partial f}{\partial u_1}+\ldots +g_n u_n\frac{\partial f}{\partial u_n},
\end{equation}
which is strict for the following choice of the filtration in $\mathcal{A}^n$: the Newton weight of $(0,\ldots , g_i, \ldots ,0)$ in is equal to $\nu (g_i) -1$, see \cite[Theorem 4.1]{K}.
\end{proof}

\begin{remark} 
\label{rem:Omegaf}
Proposition \ref{prop:Mir} provides an isomorphism of {\em rings}, and this really depends on the special form of $f$, from which we also get $\omega :=\varphi^{-1} ([f])\in H^{2}(X_{0})\subset H_{\orb}^{2} (\mathcal{X})$ where $X_{0}$ denotes the untwisted sector and $[f] : A^{*}\longrightarrow A^{*}$ 
denotes the nilpotent endomorphism induced by the multiplication by $f$ on $\mathcal{A}_f$.
\end{remark}

Last, we connect Sabbah's mixed Hodge structures and graded Milnor rings. The following result is based on the description of Sabbah's Hodge filtration with the help of the Brieskorn lattice and a Kashiwara-Malgrange $V$-filtration and the identification of the latter with a suitably normalized Newton filtration given by \cite[Theorem 4.5]{DoSa1}.

\begin{proposition} 
[\cite{DConfin}, \cite{DoSa1}, \cite{Sab}, \cite{Sab2}] 
\label{prop:HandA}
Let $f$ be a convenient and nondegenerate Laurent polynomial and let $F^{\bullet}$ be the Hodge filtration of its Sabbah's mixed Hodge structure. Then, for $p\in\zit$ we have an isomorphism,
\begin{equation}\nonumber 
\gr_F^p H \stackrel{\sim}{\longrightarrow}  \bigoplus_{n-1-p <\alpha \leq n-p} A^{\alpha}
\end{equation}
In particular, the graded vector space $\gr_F H =\oplus_p \gr_F^p H$ is isomorphic to its graded Milnor ring $A^*$. \qed
\end{proposition}

As a by-product, we get the following (well-known) result, which can be used in order to describe the nilpotent endomorphism $N$:

\begin{proposition}
[\cite{Sab}, \cite{ScSt}, \cite{Var}]
\label{prop:Nandf} 
The maps $[f]$ and $N$ have the same Jordan normal form.
\end{proposition}
\begin{proof} We follow \cite[Theorem 7.1]{ScSt}. For $\beta\in ]0,1[$, we have the commutative diagram 
$$\begin{tikzcd}
A^{\beta +n-1-p} \arrow[rightarrow] {r}{[f]} \arrow{d}{\tau^{(n-1-p)}} &  A^{\beta +n-p}\arrow{d}{\tau^{n-p}}\\
\gr_F^p H_{\beta} \arrow[rightarrow] {r}{[N]} & \gr_F^{p-1} H_{\beta} 
\end{tikzcd}$$
thanks to Proposition \ref{prop:HandA} (we use here the fact that $[N]=-\tau\partial_{\tau} \mod F^p$ and $\partial_{\tau} [\omega ]=[-f\omega]$), hence $[f]$ and $[N]$ have the same Jordan normal form.
 By Remark \ref{rem:NGrF}, $N$ and $[N]$ have the same Jordan normal form and the assertion follows. Analogous argument for $H_0$.
\end{proof}

\subsection{R\'esum\'e}

\label{sec:Resume}

Let $P$ be a simplicial polytope containing the origin as an interior point and 
let $f$ be the Laurent polynomial defined by $f (u) =\sum_{b\in \mathcal{V}(P)} u^b$. Then $f$ is a nondegenerate and convenient and we have the isomorphisms 
\begin{equation}\label{eq:IsosPremiers}
F^p H \stackrel{\sim}{\longrightarrow} \bigoplus_{i\geq p} \gr_F^i H \stackrel{\sim}{\longrightarrow} \bigoplus_{\alpha \leq n-p} A^{\alpha}\stackrel{\sim}{\longrightarrow}  \bigoplus_{\alpha \leq n-p} H^{2\alpha}_{\orb} (\mathcal{X})
\end{equation}
and 
\begin{equation}\label{eq:Isos}
H\stackrel{\sim}{\longrightarrow} \gr_F H\stackrel{\sim}{\longrightarrow} A^*  \stackrel{\sim}{\longrightarrow} H^{2*}_{\orb} (\mathcal{X}).
\end{equation}
Under these isomorphisms, $N$, $[N]$ and $[f]$ correspond respectively to $[N]$, $[f]$ and $L_{\omega}$, the latter being the nilpotent endomorphism of 
$H_{\orb}^{2*} (\mathcal{X})$ defined by
\begin{equation}\label{eq:N}
L_{\omega} (a):=\omega\cup_{\orb} a
\end{equation}
where $\omega :=\varphi^{-1} ([f])\in H^2 (X_0)$ and $\cup_{\orb}$ denotes the orbifold cup-product.
We will denote by 
\begin{equation}\label{eq:IsoHH}
\psi :H\stackrel{\sim}{\longrightarrow} H^{2*}_{\orb} (\mathcal{X})
\end{equation}
the composition of these isomorphisms.

\begin{remark} 
\label{rem:ComputeL}
It is possible to calculate $(\ref{eq:N})$ using \cite[Proposition 3.2]{Fe} and Remark 
\ref{rem:Omegaf}. The cohomology $H^* (X_{\ell} ,\cit )$ of a twisted sector is a $H^* (X_{0} ,\cit )$-module under orbifold cup-product (therefore $L_{\omega}$ preserves each sector) and this module structure is basically the standard cup-product on each sector:  we have
$L_{\omega} (a)= \omega \wedge a\in H^{p+2} (X_{\ell})$ if $a\in H^{p} (X_{\ell})$ (this is a slight abuse of notation, $\omega$ stands for its natural pullback in $H^2 (X_{\ell})$). 
\end{remark}

\begin{remark} 
\label{rem:SpectrumMonodromyInfinity}
For a calculation of the spectrum at  infinity $\Spec_f (z) =\sum_{\alpha\in\qit} \dim A^{\alpha} z^{\alpha}$ of a Laurent polynomial $f$, one can use the isomorphism (\ref{eq:Mir}) and the description of the orbifold cohomology given in Section \ref{sec:BCSsetting}. Similarly, the isomorphisms (\ref{eq:Isos}) give information about the Jordan structure of $N$ and $[f]$.
\end{remark}

\begin{example} \label{ex:Simplices} (Weighted projective spaces) 
We will say that $P$ is a simplex in $\rit^n$ if $P$ is the convex hull of $n+1$ lattice points $v_0 ,\cdots ,v_n $ and if it contains the origin as an interior point. 
The weight of the simplex $P$ is
the tuple $(q_0 ,\cdots , q_n )$ where 
\begin{equation}\nonumber
q_i := |\det (v_0 ,\cdots , \widehat{v_i},\cdots , v_n )|
\end{equation}
for $i=0,\cdots ,n$ (we will always assume that the $q_i$'s are arranged in increasing order).
 The simplex $P$ is reduced if $\gcd (q_0 ,\cdots ,q_n )=1$ (this condition is already used in \cite[p. 2]{DoSa2}). In this case, the corresponding orbifold $\mathcal{X}_P$ is a weighted projective space.
\begin{itemize}
\item Using the isomorphisms (\ref{eq:Isos}) and the description of the corresponding Sabbah's mixed Hodge structure given in \cite{DoSa2}, it is possible to give informations about $H^{2*}_{\orb} (\mathcal{X}_P)$ 
 (see also \cite{CCLT}, \cite{DoMa}, \cite{Mann1}...)
Let   
$$F:=\left\{\frac{j}{q_{i}}|\, 0\leq j\leq q_{i}-1,\ 0\leq i\leq n\right\}$$
and let $f_{1},\cdots , f_{k}$ be the elements of $F$ arranged by increasing order.
Define      
\begin{align}\nonumber
  S_{f_i}:=\{j|\ q_{j}f_i \in\zit\}\subset \{0,\cdots ,n\}\ \mbox{and}\
d_{i}:=\card S_{f_{i}}.
\end{align}
Then, the sectors are labelled by $F$ and the age of the sector $X_{f_{i}}$ for $f_{i}\in F$ is
\begin{equation}\label{eq:AgeSector}
i_{f_{i}}=d_1 +\cdots +d_{i-1}-\mu f_{i}\ \mbox{and}\ n_{f_i}:=\dim X_{f_i}= d_{i}-1
\end{equation}
for $i=2,\cdots ,k$, where $\mu := q_0 +\ldots +q_n$. 
\item Each sector $X_{\ell}$ produces exactly one Jordan block of size $n_{\ell}+1$ of the nilpotent endomorphism $N$ of $H$ since, by \cite{Kawa}, $\dim H^{2i} (X_{\ell}) =1$ and $H^{2i+1}(X_{\ell}) =0$ for all $\ell\in F$. Moreover, the spectrum at infinity is equal the disjoint union over $\ell\in F$ of the sets $\{i_{\ell}, i_{\ell}+1, \ldots , i_{\ell}+n_{\ell}\}$.
This gives another interpretation of the results obtained in \cite{DoSa2}. 
\end{itemize}
It follows that the ages of the sectors are integers if and only if $P$ is a reflexive polytope, see \cite[Proposition 5.1]{Conrads}. 

As an application, we consider the following simplex $\Delta$ which is studied by Payne in \cite[Example 4.2]{Payne} from a combinatorial perspective. We handle it using the previous description. Let $b,k$ (resp. $r$) be positive (resp. nonnegative) integers and let $\Delta$ be the reduced simplex of weight
$$(1,\ldots ,1, b,\ldots ,b)$$
where $1$ is counted $bk$-times and $b$ is counted $r+1$-times in $\rit^n$, where $n=bk+r$. We have
$$\mu =bk +b(r+1)= b(k+r+1)$$
and, since $b$ divides $\mu$, the simplex $\Delta$ is reflexive. From now on, we assume that $b\geq 2$. Then,
\begin{itemize}
\item the twisted sectors are labelled by $F=\{ 1/b, 2/b, \ldots, (b-1)/b\}$, 
\item the age of the sector $i/b$ for $1\leq i\leq b-1$ is
$$n+1+(i-1)(r+1)-\mu i/b = n+1+(i-1)(r+1)-i(k+r+1)=k(b-i),$$
\item the conjugate of the sector $\ell =i/b$ is $\ell^{-1} =(b-i)/b$ and we have 
$$i_{\ell} +i_{\ell^{-1}} =kb =n-r$$
and this is consistent with the equality $n_{\ell} =r$.
\end{itemize}
It follows that, for $f$ defined as in (\ref{eq:fOrbi}),
\begin{equation}\nonumber 
\Spec_f (z)=1+z+z^2 +\ldots +z^n +(1+\ldots +z^r) (z^{k(b-1)} +z^{k(b-2)}+\ldots +z^{k})
\end{equation}
where the product on the right hand term is the contribution of the twisted sectors: the contribution of the sector $i/b$ is 
$$k(b-i), k(b-i)+1 ,\ldots, k(b-i)+r.$$ 
This formula matches with the formula for the $h^*$-vector found in \cite[Example 4.2]{Payne} (see \cite[Theorem 6.1]{DConfin}). 
Taking $r=0$, $b\geq 3$ and $k\geq 2$, we see that the spectrum at infinity is not always unimodal. Last, and by \cite[Section 3]{DoSa2}, the nilpotent endomorphism $N$ has $b-1$ Jordan blocks of size $r+1$ and one Jordan block of size $n+1$. 
\end{example}

\section{Mixed Hodge structures on orbifold cohomology}

\label{sec:MixedHodgeOrbi}

Let $P$ be a simplicial lattice polytope containing the origin as an interior point and let $f$ be the convenient and nondegenerate Laurent polynomial defined on $(\cit^*)^n$ by $f (u) =\sum_{b\in \mathcal{V}(P)} u^b$ as in (\ref{eq:fOrbi}). In this section, we construct mixed Hodge structures on 
$H^{2*}_{\orb}(\mathcal{X})$ which correspond to the mixed Hodge structures $\MHS_f$ via the isomorphism (\ref{eq:IsoHH}).

\subsection{Filtrations}

\label{sec:Filtrations}

The definition of the following Hodge and weight filtrations
was originally suggested 
by the description given in \cite{DoSa2} of the Sabbah's mixed Hodge structures associated with 
Laurent polynomials whose Newton polytopes are reduced simplices. 
 In a different combinatorial framework, analogous filtrations are defined by Wang in \cite{Wa}. 

\subsubsection{The Hodge filtration}

\label{sec:HodgeFiltration}

We define, for $\ell\in F$ and $p\in\zit$,  
\begin{equation}\label{eq:HodgeOrbifold1}
F^p H^{2*}(X_{\ell}) :=\oplus_{i\geq p} H^{2(n-i-i_{\ell})}(X_{\ell})\ \mbox{if}\ i_{\ell}\in\zit ,
\end{equation}
and
\begin{equation}\label{eq:HodgeOrbifold2}
F^p H^{2*}(X_{\ell}) :=\oplus_{i\geq p} H^{2(n-1-i-[i_{\ell}])}(X_{\ell})\ \mbox{if}\ i_{\ell}\notin\zit.
\end{equation}
We define the following Hodge filtration on orbifold cohomology:

\begin{definition}\label{def:HodgeFiltrationOrbifold}
The Hodge filtration on $H^{2*}_{\orb}(\mathcal{X})$ is the decreasing filtration defined by
\begin{equation}\label{eq:HodgeOrbifold3}
F^p H_{\orb}^{2*}(\mathcal{X}) :=\bigoplus_{\ell\in F}F^p  H^{2*}(X_{\ell})
\end{equation}
for $p\in \zit$.
\end{definition}

\noindent This Hodge filtration has a very simple description:

\begin{proposition}\label{prop:DescGradedOrbi}
Let $p\in\zit$. Then, $F^{n-p} H_{\orb}^{2*} (\mathcal{X})=\bigoplus_{\alpha \leq p}  H_{\orb}^{2\alpha }(\mathcal{X})$.
\end{proposition}
\begin{proof}
 We have
\begin{equation}\nonumber
\bigoplus_{p-1 <\alpha \leq p}  H_{\orb}^{2\alpha } (\mathcal{X})=\oplus_{i_{\ell}\in\zit} H^{2(p-i_{\ell})} (X_{\ell})
\bigoplus\oplus_{i_{\ell}\notin\zit} H^{2(p-1- [i_{\ell}])} (X_{\ell})
\end{equation}
and the result follows.
\end{proof}

\noindent Note that $F^0 H_{\orb}^{2*}(\mathcal{X})=H_{\orb}^{2*}(\mathcal{X})$ and $F^{n+1} H_{\orb}^{2*}(\mathcal{X})=0$.

\subsubsection{The weight filtration}
\label{sec:Weight}

We define, for $\ell\in F$ and $m\in\zit$, 
\begin{equation}\label{eq:PoidsOrbifold1}
W_{m} H^{2*}(X_{\ell}) :=\bigoplus_{n+n_{\ell}-m\leq 2j}H^{2j}(X_{\ell})\  \mbox{if}\ i_{\ell}\in\zit 
\end{equation}
and
\begin{equation}\label{eq:PoidsOrbifold2}
W_{m} H^{2*}(X_{\ell}) :=\bigoplus_{n-1+n_{\ell}-m\leq 2j}H^{2j}(X_{\ell})\ \mbox{if}\ i_{\ell}\notin\zit .
\end{equation}
We define the following weight filtration on orbifold cohomology:

\begin{definition}\label{def:WeightFiltration}
The weight filtration on $H^{2*}_{\orb}(\mathcal{X})$ is the increasing filtration defined by 
\begin{equation}\nonumber \label{eq:PoidsOrbifold3}
W_m H_{\orb}^{2*}(\mathcal{X}) :=\bigoplus_{\ell\in F} W_m H^{2*}(X_{\ell})
\end{equation}
for $m\in\zit$.
\end{definition}

\noindent Recall the nilpotent endomorphism $L_{\omega}$ of $H_{\orb}^{2*} (\mathcal{X})$ defined by (\ref{eq:N}).

\begin{proposition}\label{prop:WN}
The weight filtration $W_{\bullet}$ is the weight filtration of $L_{\omega}$,
centered at $n$ on $H_{\orb}^{2*} (\mathcal{X})_0$ and centered at $n-1$ on $H_{\orb}^{2*} (\mathcal{X})_{\neq 0}$.  
\end{proposition}
\begin{proof} 
Note first that
$W_{2n} H_{\orb}^{2*}(\mathcal{X})_0 =
H_{\orb}^{2*}(\mathcal{X})_0\ \mbox{and}\ W_{-1} H_{\orb}^{2*}(\mathcal{X})_0=0$
because $n_{\ell}\leq n$ if $i_{\ell}\in\zit$ and
$W_{2(n-1)} H_{\orb}^{2*}(\mathcal{X})_{\neq 0} =
H_{\orb}^{2*}(\mathcal{X})_{\neq 0}\ \mbox{and}\ W_{-1} H_{\orb}^{2*}(\mathcal{X})_{\neq 0}=0$
because $n_{\ell}\leq n-1$ if $i_{\ell}\notin\zit$. Using the characteristic properties of the weight filtration of a nilpotent endomorphism, it is enough to check that 
\begin{itemize}
\item $L_{\omega}^{n+1}H_{\orb}^{2*} (\mathcal{X})_0=0$, $L_{\omega}(W_k H_{\orb}^{2*} (\mathcal{X})_0)\subset W_{k-2}H_{\orb}^{2*} (\mathcal{X})_0$ and that we have the isomorphisms 
\begin{equation}\label{eq:CarWeight0}
L_{\omega}^{k}:\gr^W_{n+k} H_{\orb}^{2*} (\mathcal{X})_{0} \stackrel{\sim}{\longrightarrow}\gr^W_{n-k} H_{\orb}^{2*} (\mathcal{X})_0
\end{equation}
for $k \leq n/2$,
\item $L_{\omega}^{n}H_{\orb}^{2*} (\mathcal{X})_{\neq 0}=0$, $L_{\omega}(W_k H_{\orb}^{2*} (\mathcal{X})_{\neq 0})\subset W_{k-2}H_{\orb}^{2*} (\mathcal{X})_{\neq 0}$ and that we have the isomorphisms 
\begin{equation}\label{eq:CarWeightNeq0}
L_{\omega}^{k}:\gr^W_{n-1+k}H_{\orb}^{2*} (\mathcal{X})_{\neq 0} \stackrel{\sim}{\longrightarrow}\gr^W_{n-1-k} H_{\orb}^{2*} (\mathcal{X})_{\neq 0}
\end{equation}
for $k \leq (n-1)/2$.
\end{itemize}
In both cases, the two first properties follow immediately from the definition of $L_{\omega}$.
Let us show the third on for $H_{\orb}^{2*} (\mathcal{X})_0$ (similar argument for $H_{\orb}^{2*} (\mathcal{X})_{\neq 0}$). Recall that $L_{\omega}$ preserves each sector (see Remark \ref{rem:ComputeL}).
By the hard Lefschetz property for the toric simplicial variety $X_{\ell}$ (see for instance \cite[Theorem 12.5.8]{CLS} or \cite[Theorem 2.2]{Fe} and the references therein), we have 
$$L_{\omega}^{k} : H^{n_{\ell}-k}(X_{\ell})\stackrel{\sim}{\rightarrow} H^{n_{\ell}+k}(X_{\ell})$$
and therefore
$$L_{\omega}^{k}:\gr^W_{n+k}H_{\orb}^{2*} (\mathcal{X})_{0} \stackrel{\sim}{\longrightarrow}\gr^W_{n-k}H_{\orb}^{2*} (\mathcal{X})_0$$
since $\gr^W_{n+k} H^{2*}_{\orb} (\mathcal{X})_{0} =\oplus_{i_{\ell}\in\zit} H^{n_{\ell}-k}(X_{\ell})$. 
\end{proof}

\subsubsection{Comparison with Sabbah's filtrations}

The previous filtrations match with the ones of the mixed Hodge structure $\MHS_f $:

\begin{proposition} 
\label{prop:FWSabCompare}
The Hodge filtration (resp. the weight filtration) on the orbifold cohomology defined in Definition \ref{def:HodgeFiltrationOrbifold} (resp. in Definition \ref{def:WeightFiltration}) is the image,
under the isomorphism (\ref{eq:IsoHH}),
of the Hodge filtration (resp. of the weight filtration) of the mixed Hodge structure $\MHS_f$.
\end{proposition}
\begin{proof}
By (\ref{eq:IsosPremiers}) and Proposition \ref{prop:DescGradedOrbi}, the Hodge filtration on the orbifold cohomology is the image of Sabbah's Hodge filtration by
the isomorphism (\ref{eq:IsoHH}).  
The assertion about the weight filtration follows from Proposition \ref{prop:WN}.   
\end{proof}

\subsection{Picard-Lefschetz conjugations}
\label{sec:Conjug}

Let us now define a conjugation on $H_{\orb}^{2*}(\mathcal{X})$. It turns out that the correct ones to consider are the following:

\begin{definition}
We will call  Poincar\'e-Lefschetz conjugation on $H_{\orb}^{2*}(\mathcal{X})$ any
conjugation satisfying
\begin{equation}\label{eq:ConjOrbCoh}
\overline{H^{2i} (X_{\ell})}\subset \oplus_{j\geq i} H^{2j} (X_{\ell^{-1}})
\end{equation}
for each sector $X_{\ell}$, $\ell\in F$.
\end{definition}

\begin{lemma}
Let $W_{\bullet}$ be the weight filtration on the orbifold cohomology defined in Definition \ref{def:WeightFiltration}. Then, $\overline{W_m H_{\orb}^{2*}(\mathcal{X})} =W_m H_{\orb}^{2*}(\mathcal{X})$ for any Poincar\'e-Lefschetz conjugation.
In particular,
$$W_m H_{\orb}^{2*}(\mathcal{X})=(W_m H_{\orb}^{2*}(\mathcal{X})\cap  H_{\orb ,\rit }^{2*}(\mathcal{X}))\otimes \cit$$
where $H_{\orb , \rit }^{2*}(\mathcal{X})$ denotes the real form defined by the Poincar\'e-Lefschetz conjugation involved.
\end{lemma}
\begin{proof}
Indeed, by the very definition of a Poincar\'e-Lefschetz conjugation, and since $n_{\ell} =n_{\ell^{-1}}$, we have $\overline{W_m H^{2*}(X_{\ell})} \subset W_m H^{2*}(X_{\ell^{-1}})$
for all $\ell\in F$. 
\end{proof}

 \noindent On the way, we will consider the following Poincar\'e-Lefschetz conjugations:
\begin{itemize}
\item the Poincar\'e-Lefschetz conjugation defined by
\begin{equation}\label{eq:BasicConjugForm}
\overline{a}:=(-1)^i I^* ( c (a))\in H^{2i} (X_{\ell^{-1}})
\end{equation}
for $a\in H^{2i} (X_{\ell},\cit)$,
\item the Poincar\'e-Lefschetz conjugation defined by 
\begin{equation}\label{eq:BasicConjug}
\overline{a}:= I^* ( c (a))\in H^{2i} (X_{\ell^{-1}})
\end{equation}
for $a\in H^{2i} (X_{\ell},\cit)$.
\end{itemize}
Here, $c$ denotes the standard conjugation on $H^{*} (X_{\ell})=H^{*} (X_{\ell};\rit)\otimes\cit$ 
and $I$ is the complexification of the involution 
$I:H^{*} (X_{\ell^{-1}}; \rit)\rightarrow H^{*} (X_{\ell}; \rit)$
induced by the involution $I$ considered in Section \ref{sec:BCSsetting}. The following lemma motivates the definition of the conjugations (\ref{eq:BasicConjugForm})
and (\ref{eq:BasicConjug}):

\begin{lemma} \label{lemma:Nreal}
Assume that $H_{\orb}^{2*} (\mathcal{X})$ is equipped with the conjugation (\ref{eq:BasicConjugForm}) (resp. (\ref{eq:BasicConjug})). Then, 
\begin{enumerate}
\item $I^* (L_{\omega} (a))=L_{\omega} (I^* (a))$,
\item $\overline{L_{\omega} (a)}=-L_{\omega} (\overline{a})$ (resp. $\overline{L_{\omega} (a)}=L_{\omega} (\overline{a})$),
\item $\overline{2i\pi L_{\omega} (a)}=2i\pi L_{\omega} (\overline{a})$ (resp. $\overline{2i\pi L_{\omega} (a)}=-2i\pi  L_{\omega}  (\overline{a})$).
\end{enumerate}
\end{lemma}
\begin{proof} The first equality follows from Remark \ref{rem:ComputeL}. 
We consider the remaining statements when $H_{\orb}^{2*} (\mathcal{X})$ is equipped with the conjugation (\ref{eq:BasicConjugForm}).
Since $L_{\omega}$ is real with respect to the standard conjugation $c$,
we get
 $$I^* (\overline{L_{\omega} (a)})=(-1)^{i+1}c(L_{\omega} (a))=(-1)^{i+1}L_{\omega} (c(a))=-L_{\omega} (I^* (\overline{a}))=-I^* (L_{\omega} (\overline{a})).$$ 
This shows the second point and the third assertion follows.
\end{proof}

\noindent Note that the conjugation (\ref{eq:BasicConjug}) does not match with the conjugations appearing in Sabbah's mixed Hodge structures, for which $2i\pi L_{\omega}$ is real.

\subsection{Application: construction of mixed Hodge structures}

\label{sec:MHS}

We are now ready to state:

\begin{theorem}
\label{theo:MHSOrbifold}
Assume that  $H_{\orb}^{2*} (\mathcal{X})$ is equipped with a Poincar\'e-Lefschetz conjugation and let $H_{\orb ,\rit}^{2*}(\mathcal{X})$ be the corresponding real form. 
 Then, the tuple
\begin{equation}\nonumber 
\MHS_{\orb}:=( H_{\orb}^{2*} (\mathcal{X}),  H_{\orb ,\rit}^{2*}(\mathcal{X}),  F^{\bullet} ,  W_{\bullet})
\end{equation}
is a mixed Hodge structure\footnote{ Strictly speaking, $\MHS_{\orb}$ should be indexed by the Poincar\'e-Lefschetz conjugation involved, but there will be no risk of confusion in what follows.}. \qed 
\end{theorem}

\noindent Recall the decomposition $H_{\orb}^{2*} (\mathcal{X})=H_{\orb}^{2*} (\mathcal{X})_0 \oplus H_{\orb}^{2*}(\mathcal{X})_{\neq 0}$.
  We will denote by $H_{\orb, \rit}^{2*}(\mathcal{X})_0$ and $H_{\orb, \rit}^{2*}(\mathcal{X})_{\neq 0}$ the corresponding real forms.
We construct a mixed Hodge structure on $H_{\orb}^{2* } (\mathcal{X} )_0$ (resp. 
$H_{\orb}^{2* } (\mathcal{X})_{\neq 0}$) in Proposition \ref{prop:IDegalIDel} (resp. Proposition \ref{prop:IDegalIDelAntiRef}). The theorem will follow.\\

\noindent{\em Bigradings.} By \cite[Theorem 2.13]{CKS}, a mixed Hodge structure $(H_{\cit}, H, F^{\bullet}, W_{\bullet})$ is a bigrading
\begin{equation}\label{eq:Bigrading}
H_{\cit} =\bigoplus_{i,j\in\zit} I^{i,j}
\end{equation}
such that 
\begin{equation}\label{eq:ConjugIpqMHS}
I^{i,j}=\overline{I^{j,i}}\ \mbox{mod.}\ \oplus_{a<i, b<j} I^{a,b}.
\end{equation}
The Hodge filtration and the weight filtration satisfy
\begin{equation}\label{eq:FiltrationsBigrad}
F^{p} =\bigoplus_{i\geq p} I^{i,j}\ \mbox{and}\ \ W_{m} =\bigoplus_{i+j\leq m} I^{i,j}.
\end{equation}
Given a mixed Hodge structure $(H_{\cit}, H, F^{\bullet}, W_{\bullet})$, we have 
\begin{equation}\label{def:DeligneIpq}
I^{p,q}= F^p \cap W_{p+q} \cap ( \overline{F^q}\cap W_{p+q}+\sum_{j\geq 0}  \overline{F^{q-j-1}}\cap W_{p+q-2-j}
\end{equation}
if (\ref{eq:ConjugIpqMHS}) is satisfied, which are the subspaces constructed by Deligne in \cite{Del}: a bigrading satisfying (\ref{eq:ConjugIpqMHS}) is unique and we will call it {\em the} bigrading of the mixed Hodge structure $(H_{\cit}, H, F^{\bullet}, W_{\bullet})$. A mixed Hodge structure is said to {\em split over} $\rit$ if it admits a bigrading $\oplus_{i,j\in\zit} J^{i,j}$ such that $J^{i,j}=\overline{J^{j,i}}$: if such a bigrading exists, we have $J^{i,j}=I^{i,j}= F^i \cap \overline{F^j}\cap W_{i+j}$.

\begin{example}\label{ex:BasicExample}
Let $H=H^* (X, \rit )$ be the cohomology ring of a $n$-dimensional projective manifold $X$. The bigrading is defined by $I^{p,q} =H^{n-q, n-p} (X)$ and we have
$W_m =\oplus_{k\geq 2n-m}H^k (X; \cit )$ and  $F^p =\oplus_r \oplus_{s\leq n-p}H^{r,s}(X)$.
The corresponding mixed Hodge structure splits over $\rit$.
\end{example}

\noindent {\em The reflexive case.}
We first work with $H_{\orb}^{2*} (\mathcal{X})_0$: we have $i_{\ell}\in\zit$ in the formulas below. Let us define 
\begin{equation}\label{eq:BigradMHSOrb}
I^{i,j}_{\orb ,0}:= \bigoplus_{\ell\in F,\ i_{\ell}\in\zit,\ i_{\ell}-i_{\ell^{-1}}=j-i} H^{2(n-i_{\ell}-i)}(X_{\ell})
\end{equation}
for $i,j\in\zit$ (recall that $H^* (X_{\ell})$ stands for $H^* (X_{\ell}; \cit)$; to simplify the notations, we will also write $I^{i,j}_{\orb,0}= \bigoplus\limits_{i_{\ell}-i_{\ell^{-1}}=j-i} H^{2(n-i_{\ell}-i)}(X_{\ell})$).

\begin{proposition}\label{prop:IDegalIDel}
The subspaces $I^{i,j}_{\orb,0}$ define a bigrading on $H_{\orb}^{2*}(\mathcal{X})_0$ such that 
\begin{enumerate}
\item $F^p H_{\orb }^{2*}(\mathcal{X})_0=\bigoplus\limits_{i\geq p} I^{i,j}_{\orb,0}$,
\item $W_m H_{\orb }^{2*}(\mathcal{X})_0=\bigoplus\limits_{i+j\leq m} I^{i,j}_{\orb,0}$,
\item $I^{j,i}_{\orb,0}=\overline{I^{i,j}_{\orb,0}}$ $\mod. \bigoplus\limits_{a<j, b<i}I^{a,b}_{\orb,0}$.
\end{enumerate}
In particular, the tuple 
$$\MHS_{\orb, 0}:=(H_{\orb}^{2*}(\mathcal{X})_0, H_{\orb, \rit}^{2*}(\mathcal{X})_0,  F^{\bullet} ,  W_{\bullet})$$
 is a mixed Hodge structure and $I^{i,j}_{\orb,0}=I^{i,j}$ where $I^{i,j}$ is the subspace defined by (\ref{def:DeligneIpq}).
\end{proposition}
\begin{proof}
We have
\begin{equation}\nonumber 
\bigoplus_{i\geq p} I^{i,j}_{\orb,0}=\bigoplus_{i\geq p} \bigoplus_{i_{\ell}-i_{\ell^{-1}}=j-i} H^{2(n-i_{\ell}-i)}(X_{\ell})=\bigoplus_{i_{\ell}\in\zit}\bigoplus_{i\geq p} H^{2(n-i_{\ell}-i)}(X_{\ell})=F^p H_{\orb }^{2*}(\mathcal{X})_0
\end{equation}
whence the first statement. For the weight filtration, we have
\begin{equation}\nonumber 
W_m H_{\orb }^{2*}(\mathcal{X})_0 =\bigoplus_{i_{\ell}\in\zit}\bigoplus_{n+n_{\ell}-m\leq 2k}H^{2k}(X_{\ell})
=\bigoplus_{i_{\ell}\in\zit}\bigoplus_{n+n_{\ell}-m\leq 2(n-i_{\ell}-i)}H^{2(n-i_{\ell}-i)}(X_{\ell})
\end{equation}
\begin{equation}\nonumber
=\bigoplus_{i_{\ell}\in\zit}\bigoplus_{i_{\ell}-i_{\ell^{-1}}+2i\leq m}H^{2(n-i_{\ell}-i)}(X_{\ell})
=\bigoplus_{i_{\ell}\in\zit}\bigoplus_{i+j\leq m} \bigoplus_{i_{\ell}-i_{\ell^{-1}}=j-i}H^{2(n-i_{\ell}-i)}(X_{\ell})
=\bigoplus_{i+j\leq m} I^{i,j}_{\orb,0}
\end{equation}
because $n+n_{\ell}-m=2n-i_{\ell}-i_{\ell^{-1}}-m$. The third assertion follows from the definition of the Poincar\'e-Lefschetz conjugations: indeed, for $k\geq 0$ we have
\begin{equation}\label{eq:SectorConjugIpq}
\bigoplus_{i_{\ell}-i_{\ell^{-1}}=j-i} H^{2(n-i_{\ell}-i+k)}(X_{\ell^{-1}})=
\bigoplus_{i_{\ell}-i_{\ell^{-1}}=j-k-(i-k)} H^{2(n-i_{\ell}-i+k)}(X_{\ell^{-1}})
\end{equation}
\begin{equation}\nonumber 
=\bigoplus_{i_{\ell}-i_{\ell^{-1}}=j-k-(i-k)} H^{2(n-i_{\ell^{-1}}-j+k)}(X_{\ell^{-1}})
=\bigoplus_{i_{\ell^{-1}}-i_{\ell}=j-k-(i-k)} H^{2(n-i_{\ell}-j+k)}(X_{\ell})
\subset I^{j-k,i-k}_{\orb,0}
\end{equation}
thus $I^{j,i}_{\orb,0}+\bigoplus\limits_{a<j, b<i}I^{a,b}_{\orb,0} \supset \overline{I^{i,j}_{\orb,0}}+\bigoplus\limits_{a<j, b<i}I^{a,b}_{\orb,0}$, and the reverse inclusion follows by conjugation. 
Last, and since the bigrading $I^{i,j}_{\orb,0}$ satisfies (\ref{eq:ConjugIpqMHS}) by the third point, the remaining statement follows from \cite[Theorem 2.13]{CKS}.
\end{proof}

\begin{corollary}\label{coro:MHSOrbSplit}
The mixed Hodge structure $\MHS_{\orb, 0}$ splits over $\rit$ if and only if the Poincar\'e-Lefschetz conjugation involved satisfies
$$\overline{H^{2i} (X_{\ell})}=H^{2i} (X_{\ell^{-1}})$$
 for all $\ell\in F$ and all $i\geq 0$.
\end{corollary}
\begin{proof}
Let us assume that $\overline{H^{2i}(X_{\ell})}=H^{2i}(X_{\ell^{-1}})$ for all $\ell\in F$ and all $i\geq 0$. 
We have
\begin{equation}\nonumber 
\overline{I^{p,q}}= \bigoplus_{i_{\ell}-i_{\ell^{-1}}=q-p} \overline{H^{2(n-i_{\ell}-p)}(X_{\ell})}= \bigoplus_{i_{\ell}-i_{\ell^{-1}}=q-p} H^{2(n-i_{\ell}-p)}(X_{\ell^{-1}})
\end{equation}
\begin{equation}\nonumber
= \bigoplus_{i_{\ell}-i_{\ell^{-1}}=q-p} H^{2(n-i_{\ell^{-1}}-q)}(X_{\ell^{-1}})
= \bigoplus_{i_{\ell^{-1}}-i_{\ell}=q-p} H^{2(n-i_{\ell}-q)}(X_{\ell})
=I^{q,p}
\end{equation}
and the mixed Hodge structure $\MHS_{\orb, 0}$ splits over $\rit$.
Conversely, if $j\geq 1$ we have by (\ref{eq:SectorConjugIpq})
$$\bigoplus_{i_{\ell}-i_{\ell^{-1}}=q-p} H^{2(n-i_{\ell}-p+j)}(X_{\ell^{-1}})= \bigoplus_{i_{\ell}-i_{\ell^{-1}}=p-q} H^{2(n-i_{\ell}-q+j)}(X_{\ell})\subset I^{q-j,p-j}$$
therefore if $\overline{H^{2i} (X_{\ell})}\neq H^{2i} (X_{\ell^{-1}})$ for some $\ell\in F$ the mixed Hodge structure $\MHS_{\orb ,0 }$ does not split over $\rit$.
\end{proof}

\noindent {\em The anti-reflexive case.}
We now work on $H_{\orb}^{2*} (\mathcal{X})_{\neq 0}$. Let us define, for $i,j\in\zit$,
\begin{equation}\label{eq:BigradMHSOrbNeq0}
I^{i,j}_{\orb ,\neq 0}:= \bigoplus_{\ell\in F,\ i_{\ell}\notin\zit ,\ [i_{\ell}]-[i_{\ell^{-1}}]=j-i} H^{2(n-1-[i_{\ell}]-i)}(X_{\ell}).
\end{equation}
We have the following counterpart of Proposition \ref{prop:IDegalIDel} (its proof is analogous):

\begin{proposition}\label{prop:IDegalIDelAntiRef}
The subspaces $I^{i,j}_{\orb, \neq 0}$ define a bigrading on $H_{\orb}^{2*}(\mathcal{X})_{\neq 0}$ such that
\begin{enumerate}
\item $F^p H_{\orb }^{2*}(\mathcal{X})_{\neq 0}=\bigoplus\limits_{i\geq p} I^{i,j}_{\orb,\neq 0}$,
\item $W_m H_{\orb}^{2*}(\mathcal{X})_{\neq 0}=\bigoplus\limits_{i+j\leq m} I^{i,j}_{\orb, \neq 0}$,
\item $I^{j,i}_{\orb, \neq 0}=\overline{I^{i,j}_{\orb,0}}$ $\mod. \bigoplus\limits_{a<j, b<i}I^{a,b}_{\orb,\neq 0}$.
\end{enumerate}
In particular, the tuple 
$$\MHS_{\orb, \neq 0}:=(H_{\orb}^{2*}(\mathcal{X})_{\neq 0}, H_{\orb, \rit}^{2*}(\mathcal{X})_{\neq 0},  F^{\bullet} ,  W_{\bullet})$$
 is a mixed Hodge structure and $I^{i,j}_{\orb, \neq 0}=I^{i,j}$
where $I^{i,j}$ is the subspace defined by (\ref{def:DeligneIpq}).\qed 
\end{proposition}

\noindent The counterpart of Corollary \ref{coro:MHSOrbSplit} is straightforward.

\subsection{Hodge numbers}
\label{sec:HodgeNumbers}

We now compute the Hodge numbers
$h^{p,q}:=\dim \gr_F^p \gr^W_{p+q} H_{\orb}^{2*}(\mathcal{X})$
of the mixed Hodge structures $\MHS_{\orb}$. According to the decomposition (\ref{eq:NotationDecompHorb}), we will write $h^{p,q}=h_{0}^{p,q}+ h_{\neq 0}^{p,q}$.

\begin{proposition} \label{prop:GraduesOrbi}
We have
\begin{equation} \nonumber 
h_{0}^{p,q} =\sum_{i_{\ell}\in\zit ,\ i_{\ell}-i_{\ell^{-1}}=q-p}
\dim H^{2(n_{\ell} +i_{\ell^{-1}}-p)}(X_{\ell})
\end{equation}
and
\begin{equation}\nonumber 
h_{\neq 0}^{p,q} =\sum_{i_{\ell}\notin\zit ,\ [i_{\ell}]-[i_{\ell^{-1}}]=q-p}
\dim H^{2(n_{\ell} +[i_{\ell^{-1}}]-p)}(X_{\ell} )
\end{equation}
for $p,q\in\zit$.
\end{proposition}
\begin{proof} The assertions follow from Proposition \ref{prop:IDegalIDel}, Proposition \ref{prop:IDegalIDelAntiRef} and equality (\ref{eq:FormuleDim}). 
\end{proof}

\begin{remark} We have $h^{0,q}_0 =0$ if $q>0$.
\end{remark}

Recall that the mixed Hodge structure $\MHS_{\orb}$ is of Hodge-Tate type if 
$h^{p,q}=0$ for $p\neq q$ (see also Section \ref{sec:HThpq} below). The following statement can be found in \cite{D13} with a different proof:

\begin{corollary}\label{coro:hpp}
The mixed Hodge structure $\MHS_{\orb}$ is of Hodge-Tate type if and only if 
$[i_{\ell}]=[i_{\ell^{-1}}]\ \mbox{ for all}\ \ell\in F$. 
\end{corollary}
\begin{proof} The if part follows immediately from Proposition \ref{prop:GraduesOrbi}. Conversely, and again by Proposition \ref{prop:GraduesOrbi},
we have $h_{\neq 0}^{[i_{\ell^{-1}}],[i_{\ell}]} \geq \dim H^{2n_{\ell}}(X_{\ell} )>0$ 
and $h_{0}^{i_{\ell^{-1}}, i_{\ell}} \geq \dim H^{2n_{\ell}}(X_{\ell} )>0$
therefore  if there exists a label $\ell$ such that $[i_{\ell}]\neq [i_{\ell^{-1}}]$ the mixed Hodge structure
$\MHS_{\orb}$ is not of Hodge-Tate type.
\end{proof}

\begin{remark}
By Proposition \ref{prop:FWSabCompare}, the Hodge numbers 
of the mixed Hodge structures $\MHS_f$ and $\MHS_{\orb}$ are the same.
 A formula for the Hodge numbers of the mixed Hodge structure $\MHS_f$ is also given in \cite[Corollary 4.10]{Wa} in terms of the $h$-vector of the Newton polytope $P$ of $f$.
\end{remark}

\begin{remark} \label{rem:PropertiesHN}
Using Proposition \ref{prop:GraduesOrbi}, we easily get expected properties for the Hodge numbers of the mixed Hodge structure $\MHS_{\orb}$ and $\MHS_f$ (see for instance \cite[Section 2]{DConfin}). For instance,
\begin{enumerate}
\item  $h^{p,q}_0=h^{q,p}_0$ and $h^{p,q}_{\neq 0}=h^{q,p}_{\neq 0}$,
\item $h^{n-p,n-q}_0=h^{p,q}_0$ and $h^{n-1-p,n-1-q}_{\neq 0}=h^{p,q}_{\neq 0}$.
\end{enumerate}
The first point follows  from the equality $\dim H^{2k}(X_{\ell})=\dim H^{2k}(X_{\ell^{-1}})$ and the second one from the Poincar\'e duality for the simplicial toric varieties $X_{\ell}$.
\end{remark}

\begin{remark}
The previous Hodge numbers could be compared with the irregular Hodge numbers defined in \cite{SY} and their orbifold version considered in \cite{HL}. In the setting of mirror symmetry, the computation of Hodge numbers of Calabi-Yau varieties produced by crepant resolutions of toric varieties coming from reflexive polytopes is carried out in \cite{GKR}. 
\end{remark}

When $P$ is a reduced simplex, we get a formula depending only on the ages of the sectors (recall that these ages are easily calculated from the weight of the simplex, see Example \ref{ex:Simplices}).

\begin{corollary}
\label{coro:NbreHodgeSimplex}
Let $P$ be a reduced simplex. Then,
\begin{equation}\nonumber
h_{\neq 0}^{p,q} =\Card \{ \ell\in F, \ [i_{\ell}]-[i_{\ell^{-1}}]=q-p,\ [i_{\ell}]\leq n-1-p,\ [i_{\ell^{-1}}]\leq p \}
\end{equation}
and
\begin{equation}\nonumber 
h_{0}^{p,q} =\Card \{ \ell\in F, \ i_{\ell}-i_{\ell^{-1}}=q-p,\ i_{\ell}\leq n-p,\ i_{\ell^{-1}}\leq p\}
\end{equation}
for $p,q\in\zit$.
\end{corollary}
\begin{proof}
By \cite{Kawa}, we have $\dim H^{2i}(X_{\ell} )=1$ and 
$\dim H^{2i+1}(X_{\ell} )=0$ for $i=0,\ldots , n_{\ell}$ and we get the formulas from Proposition \ref{prop:GraduesOrbi}. The restrictions about the ages of the sectors follow from the inequalities $n_{\ell}\geq n_{\ell} +[i_{\ell^{-1}}]-p\geq 0$ and (\ref{eq:FormuleDim}).
\end{proof}

\subsection{Examples}
We give three examples: the first (resp. the third) one deals with the Hodge-Tate (resp. the reflexive) case 
and the second one deals with the general case.

\begin{example} \label{ex:HodgeNumbersn5}
Let $f$
be the Laurent polynomial defined on $(\cit^*)^5$ by
$$f (u_1 , u_2 ,u_3 , u_4 ,u_5 )=u_1+u_2 +u_3 +u_4 + u_5 +\frac{1}{u_1^2 u_2^2 u_3^3 u_4^3 u_5^3}.$$
Its Milnor number is $\mu =14$. Its Newton polytope $P$ is the reduced simplex of weight $(1,2,2,3,3,3)$ and $H^{2*}_{\orb} (\mathcal{X})$ is the orbifold cohomology of the weighted projective space $\ppit (1,2,2,3,3,3)$ (this example is borrowed from \cite{Jia}).
The sectors are labelled by the set
$$F=\{0, 1/3 , 1/2 , 2/3 \}$$
 with respective ages 
$i_{0}=0$, $i_{1/3}=4/3$, $i_{2/3}=5/3$, $i_{1/2}=2$ and dimensions
 $n_0 =5$, $n_{1/3} =n_{2/3} =2$ and $n_{1/2}=1$. Since $(1/3)^{-1}=2/3$ and $(1/2)^{-1}=1/2$, the mixed Hodge structure $\MHS_{\orb}$ is of Hodge-Tate type and we get the Hodge diamond
$$\begin{array}{ccccccccccc}
 & & & &   & 1 & &  &  & &  \\
& & & & 0  &  & 0 &  &  & &  \\
& & & 0 &   & 3 &  & 0  &  & &  \\
& & 0 &  & 0  &  & 0 &   & 0  & &  \\
& 0 &  & 0 &   & 4 &  & 0  &   & 0 &  \\
0 &  & 0 &  & 0  &  & 0 &   & 0  &  & 0 \\
  & 0 &  & 0 &   & 4  &  & 0  &   & 0 &  \\
  &  & 0 &  & 0  &   & 0  &   & 0  &  &  \\
&  &  & 0 &   & 1  &   & 0  &   &  &  \\
&  &  &  & 0  &   &  0 &   &   &  &  \\
&  &  &  &   &  1 &   &   &   &  &  \\
\end{array}$$
$$=
\begin{array}{ccccccccccc}
 & & & &   & 1 & &  &  & &  \\
& & & & 0  &  & 0 &  &  & &  \\
& & & 0 &   & 1 &  & 0  &  & &  \\
& & 0 &  & 0  &  & 0 &   & 0  & &  \\
& 0 &  & 0 &   & 2 &  & 0  &   & 0 &  \\
0 &  & 0 &  & 0  &  & 0 &   & 0  &  & 0 \\
  & 0 &  & 0 &   & 2  &  & 0  &   & 0 &  \\
  &  & 0 &  & 0  &   & 0  &   & 0  &  &  \\
&  &  & 0 &   & 1  &   & 0  &   &  &  \\
&  &  &  & 0  &   &  0 &   &   &  &  \\
&  &  &  &   &  1 &   &   &   &  &  \\
\end{array}
+
\begin{array}{ccccccccccc}
 & & & &   & 0 & &  &  & &  \\
& & & & 0  &  & 0 &  &  & &  \\
& & & 0 &   & 2 &  & 0  &  & &  \\
& & 0 &  & 0  &  & 0 &   & 0  & &  \\
& 0 &  & 0 &   & 2 &  & 0  &   & 0 &  \\
0 &  & 0 &  & 0  &  & 0 &   & 0  &  & 0 \\
  & 0 &  & 0 &   & 2  &  & 0  &   & 0 &  \\
  &  & 0 &  & 0  &   & 0  &   & 0  &  &  \\
&  &  & 0 &   & 0  &   & 0  &   &  &  \\
&  &  &  & 0  &   &  0 &   &   &  &  \\
&  &  &  &   &  0 &   &   &   &  &  \\
\end{array}
$$

\noindent for the mixed Hodge structures $\MHS_{\orb}$ and $\MHS_f$, where the sum refers to the direct sum decomposition (\ref{eq:NotationDecompHorb}). 
\end{example}

\begin{example} \label{ex:HodgeNumbersn5Bis}
Let $f$
be the Laurent polynomial defined on $(\cit^*)^5$ by
$$f (u_1 , u_2 ,u_3 , u_4 ,u_5 )=u_1+u_2 +u_3 +u_4 + u_5 +\frac{1}{u_1^2 u_2^3 u_3^4 u_4^7 u_5^{11}}.$$
We have $\mu =28$. Its Newton polytope $P$ is the reduced simplex of weight $(1,2,3,4,7,11)$ and $H^{2*}_{\orb} (\mathcal{X})$ is the orbifold cohomology of the weighted projective space $\ppit (1,2,3,4,7,11)$.
The sectors are labelled by the set
$$F=\{0, 1/11, 1/7, 2/11, 1/4, 3/11,2/7, 1/3,4/11,3/7,5/11,1/2,6/11,4/7,7/11,2/3,5/7,$$
$$ 8/11, 3/4,9/11,6/7 , 10/11 \}$$
 with respective ages 
$$0, 38/11, 3, 32/11, 2, 26/11, 3, 8/3, 31/11, 2, 25/11, 2, 30/11, 3, 24/11, 7/3, 2,$$
$$ 29/11, 3, 23/11, 2 , 17/11 .$$
We get the Hodge diamond
$$\begin{array}{ccccccccccc}
 & & & &   & 1 & &  &  & &  \\
& & & & 0  &  & 0 &  &  & &  \\
& & & 0 &   & 1 &  & 0  &  & &  \\
& & 0 &  & 0  &  & 0 &   & 0  & &  \\
& 0 &  & 1 &   & 12 &  & 1  &   & 0 &  \\
0 &  & 0 &  & 4  &  & 4 &   & 0  &  & 0 \\
  & 0 &  & 0 &   & 2  &  & 0  &   & 0 &  \\
  &  & 0 &  & 0  &   & 0  &   & 0  &  &  \\
&  &  & 0 &   & 1  &   & 0  &   &  &  \\
&  &  &  & 0  &   &  0 &   &   &  &  \\
&  &  &  &   &  1 &   &   &   &  &  \\
\end{array}$$

$$=\begin{array}{ccccccccccc}
 & & & &   & 1 & &  &  & &  \\
& & & & 0  &  & 0 &  &  & &  \\
& & & 0 &   & 1 &  & 0  &  & &  \\
& & 0 &  & 0  &  & 0 &   & 0  & &  \\
& 0 &  & 0 &   & 2 &  & 0  &   & 0 &  \\
0 &  & 0 &  & 4  &  & 4 &   & 0  &  & 0 \\
  & 0 &  & 0 &   & 2  &  & 0  &   & 0 &  \\
  &  & 0 &  & 0  &   & 0  &   & 0  &  &  \\
&  &  & 0 &   & 1  &   & 0  &   &  &  \\
&  &  &  & 0  &   &  0 &   &   &  &  \\
&  &  &  &   &  1 &   &   &   &  &  \\
\end{array}
+
\begin{array}{ccccccccccc}
 & & & &   & 0 & &  &  & &  \\
& & & & 0  &  & 0 &  &  & &  \\
& & & 0 &   & 0 &  & 0  &  & &  \\
& & 0 &  & 0  &  & 0 &   & 0  & &  \\
& 0 &  & 1 &   & 10 &  & 1  &   & 0 &  \\
0 &  & 0 &  & 0  &  & 0 &   & 0  &  & 0 \\
  & 0 &  & 0 &   & 0  &  & 0  &   & 0 &  \\
  &  & 0 &  & 0  &   & 0  &   & 0  &  &  \\
&  &  & 0 &   & 0  &   & 0  &   &  &  \\
&  &  &  & 0  &   &  0 &   &   &  &  \\
&  &  &  &   &  0 &   &   &   &  &  \\
\end{array}$$

\noindent for the mixed Hodge structures $\MHS_{\orb}$ and $\MHS_f$.
\end{example}

\begin{example} \label{ex:HodgeNumbersn4}
Let $f$
be the Laurent polynomial defined on $(\cit^*)^4$ by
$$f (u)=u_1 +u_2 +u_3 +u_4 +\frac{1}{u_1^5 u_2^{12} u_3^{12} u_4^{30}}.$$
We have $\mu =60$. Its Newton polytope $P$ is the reduced and reflexive simplex of weight $(1,5,12,12,30)$ in $\rit^4$.
We get the Hodge diamond 
$$\begin{array}{ccccccccc}
 & &  & & 1 & & & &  \\
 & & & 0 & & 0 &  & & \\
 & & 0 & & 6 & & 0 & & \\
& 0 & & 5 & & 5 & & 0 &\\ 
0 & & 4 & & 18 & & 4 & & 0 \\
& 0 & & 5 & & 5 & & 0 &\\ 
 & & 0 & & 6 & & 0 & & \\
 & & & 0 & & 0 &  & & \\
 & &  & & 1 & & & &  
\end{array}$$
for the mixed Hodge structures $\MHS_{\orb}$ and $\MHS_f$. The Hodge diamond is symmetric about the horizontal axis and this agrees with Remark \ref{rem:PropertiesHN} since $P$ is reflexive.
\end{example}

\section{Interlude: the Hodge-Tate condition}

\label{sec:HThpq}

After Corollary \ref{coro:hpp}, we naturally set:

\begin{definition}
We will call the condition
\begin{equation}\label{eq:ConditionClef}
[i_{\ell}]=[i_{\ell^{-1}}]\ \mbox{ for all}\ \ell\in F 
\end{equation}
the  Hodge-Tate condition.
\end{definition}

\noindent When the ages are all integers, this condition already appears in \cite{Fe}. It also appears in this form in \cite{D13} (where the relation with a conjecture by Katzarkov, Kontsevitch and Pantev \cite{KKP} is discussed) and one can refer to loc. cit. in order to have an idea of how common this condition is.

The Hodge-Tate condition deserves attention for other various reasons. First, if it is satisfied, 
we get an easy description of the weight filtration:

\begin{proposition} 
\label{prop:iellW}
Let $W_{\bullet}$ be the weight filtration of a mixed Hodge structure $\MHS_{\orb}$ defined by Theorem \ref{theo:MHSOrbifold} and 
assume that the Hodge-Tate condition (\ref{eq:ConditionClef}) is satisfied. Then,
\begin{equation}\label{eq:HTW}
W_{2k} H_{\orb}^{2*} (\mathcal{X})=W_{2k+1}H_{\orb}^{2*} (\mathcal{X})=\oplus_{n-k-1<\alpha} H^{2\alpha}_{\orb}(\mathcal{X})
\end{equation}
for $m\in\zit$. 
 \end{proposition}
\begin{proof}
The proof of the first assertion is analogous to the one of Proposition \ref{prop:DescGradedOrbi}: one uses (\ref{eq:PoidsOrbifold1}), (\ref{eq:PoidsOrbifold2})
 and the fact that, under the assumption (\ref{eq:ConditionClef}), $n_{\ell}=n-2 i_{\ell}$ if $i_{\ell}\in\zit$ and $n_{\ell}=n-2[i_{\ell}]-1$ if $i_{\ell}\notin\zit$
in order to get
\begin{equation}\nonumber
W_m H_{\orb}^{2*} (\mathcal{X})=\oplus_{2i \leq m} \oplus_{n-i-1<\alpha\leq n-i} H_{\orb}^{2\alpha} (\mathcal{X})
\end{equation}
for $m\in\zit$.
\end{proof}

One can also use the standard conjugation $c$ defined in Section \ref{sec:Conjug} in order to get a mixed Hodge structure on $H_{\orb}^{2*} (\mathcal{X})$ (we denote by $H_{\orb}^{2*}(\mathcal{X}; \rit )$ the corresponding real form):

\begin{proposition}
\label{prop:MHSOrbifoldC}
Assume that  $H_{\orb}^{2*} (\mathcal{X})$ is equipped with the standard conjugation $c$. Assume moreover that the Hodge-Tate condition (\ref{eq:ConditionClef}) is satisfied.
 Then, the tuple 
\begin{equation}\nonumber 
\MHS_{\orb}:=( H_{\orb}^{2*} (\mathcal{X}),  H_{\orb}^{2*}(\mathcal{X}; \rit ),  F^{\bullet} ,  W_{\bullet})
\end{equation}
is a mixed Hodge structure which splits over $\rit$. 
\end{proposition}
\begin{proof} If the Hodge-Tate condition (\ref{eq:ConditionClef}) is satisfied we have $I^{i,i}_{\orb ,0}=\bigoplus_{i_{\ell}\in\zit} H^{2(n-i_{\ell}-i)}(X_{\ell})$ and $I^{i,j}_{\orb ,0}=0$ if $i\neq j$. In particular, $c(I^{i,i}_{\orb,0})=I^{i,i}_{\orb,0}$ and the statement for $H_{\orb}^{2*} (\mathcal{X})_0$ follows from the discussion in Section \ref{sec:MHS}. Similar arguments for $H_{\orb}^{2*} (\mathcal{X})_{\neq 0}$.
\end{proof}

\begin{remark} Other (orbifold) bigradings could be considered. For instance, if $H_{\orb}^{2*} (\mathcal{X})=H_{\orb}^{2*} (\mathcal{X})_0$ and after Example \ref{ex:BasicExample}, one could work with the bigrading 
\begin{equation}\nonumber
I^{i,j}_{F,0}:= H^{n-j, n-i}_{\orb, 0}:=  \bigoplus_{\ell\in F} H^{n-j-i_{\ell}, n-i-i_{\ell}}(X_{\ell})
\end{equation}
induced by the orbifold Dolbeault cohomology groups, see \cite{Fe}. Since $X_{\ell}$ is a toric simplicial variety, we have $I^{i,j}_{F,0}=0$ if $i\neq j$ and it is easily seen that $I^{i,j}_{\orb ,0}= I^{i,j}_{F,0}$ if and only if  the mixed Hodge structure $\MHS_{\orb}$ is of Hodge-Tate type. More generally, one could consider the bigrading of $H_{\orb}^{2*} (\mathcal{X})_{\neq 0}$ 
\begin{equation}\nonumber 
I^{\lambda , \mu }_{F,\neq 0}:= H^{n-\mu, n-\lambda }_{\orb}:=  \bigoplus_{\ell\in F} H^{n-\mu-i_{\ell}, n-\lambda-i_{\ell}}(X_{\ell})
\end{equation}
for $\lambda ,\mu\in\qit $, induced by the orbifold Dolbeault cohomology groups, and define the orbifold Hodge numbers $h_{\orb}^{\lambda ,\mu}:=\dim I^{\lambda , \mu }_{F,\neq 0}$. See \cite{HL} for constructions in the same spirit.
\end{remark}

\begin{remark} 
\label{rem:KnownCases}
 Assume that the Hodge-Tate condition (\ref{eq:ConditionClef}) is satisfied and that 
$H_{\orb}^{2*} (\mathcal{X})=H_{\orb}^{2*} (\mathcal{X})_0$. Then,
\begin{equation}\nonumber
F^p H_{\orb}^{2*} (\mathcal{X})=\bigoplus_{\alpha \leq n-p} H_{\orb}^{2\alpha} (\mathcal{X})
\ \mbox{and}\ W_m H_{\orb}^{2*} (\mathcal{X})=\bigoplus_{2i \leq m} H_{\orb}^{2(n-i)} (\mathcal{X}).
\end{equation}
In our (simplicial) toric setting, these filtrations match with the ones defined in \cite[Section 5]{Fe} and the mixed Hodge structure of Proposition
\ref{prop:MHSOrbifoldC} coincide with the one defined in \cite[Theorem 5.3]{Fe}.
\end{remark}

 The Hodge-Tate condition is also connected with hard Lefschetz properties for orbifold cohomology. In what follows, we write $N$ for $L_{\omega}$, the nilpotent endomorphism of 
$H_{\orb}^{2*} (\mathcal{X})$ defined by (\ref{eq:N}).
We will say that $H_{\orb}^{2*}(\mathcal{X})$ satisfies the hard Lefschetz property if the powers of $N$ induces isomorphisms
\begin{equation}\label{eq:HL0}
N^{n-2k}: H_{\orb}^{2k} (\mathcal{X}) \stackrel{\sim}{\longrightarrow} H_{\orb}^{2(n-k)} (\mathcal{X})
\end{equation}
for all integers $k$ such that $0\leq k\leq n/2$ and
\begin{equation}\label{eq:HLneq0}
N^{n-1-2k} : H^{2(\beta +k)}_{\orb }(\mathcal{X})\stackrel{\sim}{\longrightarrow}H^{2(\beta +n-1-k)}_{\orb }(\mathcal{X})
\end{equation}
for all $\beta \in ]0,1[$ and integers $k$ such that $0\leq k\leq (n-1)/2$. The next result 
 can be found with a different proof in \cite[Theorem 3.4]{D13} and is useful in order to test whether the orbifold cohomology satisfies a K\"{a}hler package (see Remark \ref{remark:Kahler Package}).

\begin{proposition}
\label{proposition:HTEquiv}
The hard Lefschetz properties (\ref{eq:HL0}) and (\ref{eq:HLneq0}) hold if and only if the Hodge-Tate condition (\ref{eq:ConditionClef}) is satisfied.
\end{proposition}
\begin{proof}
By the description of the Hodge filtration given by Proposition \ref{prop:DescGradedOrbi}, we have 
$$\gr^{n-1-k}_F H^{2*}_{\orb}(\mathcal{X})_{\neq 0}=\bigoplus_{0<\beta <1} H^{2(\beta +k)}_{\orb }(\mathcal{X})\ \mbox{and}\ \gr^{k}_F H^{2*}_{\orb}(\mathcal{X})_{\neq 0}=\bigoplus_{0<\beta <1} H^{2(\beta +n-1-k)}_{\orb }(\mathcal{X})$$
therefore (\ref{eq:HLneq0}) holds if and only if 
$$N^{n-1-2k}:\gr^{n-1-k}_F H^{2*}_{\orb}(\mathcal{X})_{\neq 0}\stackrel{\sim}{\longrightarrow}\gr^{k}_F H^{2*}_{\orb}(\mathcal{X})_{\neq 0}$$
that is if and only if $[i_{\ell}]=[i_{\ell^{-1}}]$ for all $\ell\in F$ such that $i_{\ell}\notin\zit$ by Corollary \ref{coro:hpp} and \cite[Corollary 2.3]{Sab2}.
In the same way, (\ref{eq:HL0}) is true if and only if $i_{\ell}=i_{\ell^{-1}}$ for all $\ell\in F$ such that $i_{\ell}\in\zit$.
\end{proof}

\section{Polarizations}
\label{sec:Polarization}

We now explain how the orbifold Poincar\'e duality defined by 
Chen and Ruan in \cite{CR} provides a polarization of the mixed Hodge structure $\MHS_{\orb}$ defined in
Theorem \ref{theo:MHSOrbifold}.  We consider two special cases, leading to different statements: the first one, discussed in Section \ref{sec:OrbiPoincDual}, looks like the duality of $\mathcal{D}$-modules introduced in \cite{Sab} 
(see Appendix \ref{sec:DualityDmodules})
and fits with Sabbah's mixed Hodge structures; the second one, discussed in Section \ref{sec:Alternate} and defined only in the Hodge-Tate case, is more classical and follows along the lines the approach used in \cite{Fe}.

\subsection{Orbifold Poincar\'e-Saito duality and polarization}

\label{sec:OrbiPoincDual}

In this section, we assume that $H^{2*}_{\orb} (\mathcal{X})$ is equipped with
the Poincar\'e-Lefschetz conjugation (\ref{eq:BasicConjugForm}). 
Using the orbifold Poincar\'e duality defined in \cite{CR}, we construct a polarization of the mixed Hodge structure $\MHS_{\orb}$ in the sense of Definition \ref{def:PolarizedMHS}. 
This construction is basically inspired by the description of the duality of $\mathcal{D}$-module for weighted projective spaces given in Appendix \ref{sec:DualityDmodulesWPS}. 
In what follows, we will set $r=n-1$ when considering $H^{2*}_{\orb}(\mathcal{X})_{\neq 0}$ and $r=n$ when considering $H^{2*}_{\orb}(\mathcal{X})_{0}$. 
We will denote by $H^{2*}_{\orb}(\mathcal{X})$ either $H^{2*}_{\orb}(\mathcal{X})_{0}$ or 
$H^{2*}_{\orb}(\mathcal{X})_{\neq 0}$.

\subsubsection{Orbifold Poincar\'e-Saito duality}

After \cite[Proposition 3.3.1]{CR}, and for $0\leq \alpha\leq n$, we define the bilinear form 
\begin{equation}\nonumber
g : H^{2\alpha}_{\orb} (\mathcal{X})\times H^{2(n-\alpha )}_{\orb} (\mathcal{X})\rightarrow \cit
\end{equation}
 as the direct sum over the sectors $X_{\ell}$ such that  $\alpha -i_{\ell}$ is a nonnegative integer of
\begin{equation}\nonumber 
g^{(\ell)} : H^{2(\alpha -i_{\ell})} (X_{\ell}) \times H^{2(n-\alpha-i_{\ell^{-1}})} (X_{\ell^{-1}})\rightarrow\cit 
\end{equation}
where 
\begin{equation} \label{def:OrbiBilinearDmod}
g^{(\ell)} (a,b)=\frac{(-1)^{[\alpha]}}{(2i\pi )^{r}}\int_{\mathcal{X}_{\ell}}^{\orb} a \wedge I^* (b),
\end{equation}
$I$ denoting the involution $I:X_{\ell}\rightarrow X_{\ell^{-1}}$ alluded to in Section \ref{sec:BCSsetting}.
 We refer to \cite{CR} or \cite[III 3.3]{Mann} for the definition of the orbifold integral. 
We will call $g$ the {\em Poincar\'e-Saito duality}.\\

We first list the basic properties of the Poincar\'e-Saito duality, keeping in mind polarizations of mixed Hodge structures (see Appendix \ref{sec:PolarizedMHS}).

\begin{lemma} \label{lemma:propertiesgorbifold}
Let $g$ be the Poincar\'e-Saito duality and assume that $H^{2*}_{\orb} (\mathcal{X})$ is equipped with the Poincar\'e-Lefschetz conjugation (\ref{eq:BasicConjugForm}). Then,
\begin{enumerate}
\item $g$ is nondegenerate,
\item $g$ is $ (-1)^{r}$-symmetric, 
\item $\overline{g(a,b)}= g(\overline{a}, \overline{b})$ if $a\in H^{2(\alpha -i_{\ell})} (X_{\ell})$ and $b\in H^{2(n-\alpha-i_{\ell^{-1}})} (X_{\ell^{-1}})$, 
\item $g(2i\pi Na, b)+g(a, 2i\pi Nb)=0$.
\end{enumerate}
\end{lemma}
\begin{proof}
The proof basically relies on the normalization of the integral in the definition of $g^{(\ell)}$.
\begin{enumerate} 
\item By (\ref{eq:FormuleDim}), we have $n-\alpha-i_{\ell^{-1}}=n_{\ell}-(\alpha -i_{\ell})$ therefore $g^{(\ell)}$ is isomorphic to the ordinary Poincar\'e duality on $X_{\ell}$. In particular, $g$ is nondegenerate since $g^{(\ell)}$ is. 
\item  Let $a\in H^{2(\alpha -i_{\ell})} (X_{\ell})$ and let $b\in H^{2(n-\alpha-i_{\ell^{-1}})} (X_{\ell^{-1}})$. Then,
$$g (b,a)= \frac{(-1)^{[n-\alpha]}}   {(2i\pi )^{r}}    \int_{\mathcal{X}_{\ell^{-1}}}^{\orb} b \wedge I^* (a)=\frac{(-1)^{r-[\alpha]}} {(2i\pi )^{r}} 
\int_{\mathcal{X}_{\ell}}^{\orb} I^* (b) \wedge a=  (-1)^{r} g(a,b).$$
\item Let $a\in H^{2(\alpha -i_{\ell})} (X_{\ell})$ and let $b\in H^{2(n-\alpha-i_{\ell^{-1}})} (X_{\ell^{-1}})$. Then,
$$\overline{g (a,b)}= (-1)^{r}\frac{(-1)^{[\alpha]}}   {(2i\pi )^{r}} \int_{\mathcal{X}_{\ell}}^{\orb} c(a) \wedge I^* (c(b))
$$
$$=
\frac{(-1)^{[\alpha] -[i_{\ell^{-1}}]     -[i_{\ell}]        }}   {(2i\pi )^{r}}
\int_{\mathcal{X}_{\ell^{-1}}}^{\orb} \overline{a}\wedge I^* (\overline{b})
= 
\frac{(-1)^{[\beta]}}   {(2i\pi )^{r}}
\int_{\mathcal{X}_{\ell^{-1}}}^{\orb} \overline{a}\wedge I^* (\overline{b}) =g( \overline{a} , \overline{b})$$
because $I^* (\overline{a})=(-1)^{\alpha -i_{\ell}}c(a)$, $I^* (\overline{b})=(-1)^{n-\alpha -i_{\ell^{-1}}}c(b)$, $n-i_{\ell^{-1}}-i_{\ell}=r -[i_{\ell^{-1}}]     -[i_{\ell}]$ and
$\overline{a}\in H^{2(\beta -i_{\ell^{-1}})} (X_{\ell^{-1}})$ with $\beta -i_{\ell^{-1}}
=\alpha-i_{\ell}\in\zit$.
\item Let $a\in H^{2(\alpha -i_{\ell})} (X_{\ell})$. Then,
$$g (2i\pi Na,b)=
\frac{(-1)^{[\alpha] +1}}{(2i\pi )^{r}}
\int_{\mathcal{X}_{\ell}}^{\orb} 2i\pi N(a) \wedge I^* (b)
=
\frac{(-1)^{[\alpha] +1}}{(2i\pi )^{r}}\int_{\mathcal{X}_{\ell}}^{\orb} a \wedge I^* (2i\pi N(b))$$
$$=-g(a, 2i\pi Nb)$$
where the second equality follows from Lemma \ref{lemma:Nreal}.
\end{enumerate}
This completes the proof.
\end{proof}

\begin{lemma}\label{lemma:gF}
Let $F^{\bullet}$ and $W_{\bullet}$ be the Hodge filtration of the mixed Hodge structure $\MHS_{\orb}$ and let $g$ be the Poincar\'e-Saito duality. Then,
\begin{enumerate}
\item $g(F^p , F^{r+1-p})=0$,
\item $g(W_i , W_{2r-i-1} )=0$.
\end{enumerate}
\end{lemma}
\begin{proof}
The first point follows from Proposition \ref{prop:DescGradedOrbi} and from the definition of $g$.
 For the second point, we have 
$$W_i =\oplus_{r+n_{\ell}-i\leq 2j} H^{2j}(X_{\ell})\ \mbox{and}\ 
W_{2r-i-1} =\oplus_{n_{\ell}+i+1-r\leq 2j'} H^{2j'}(X_{\ell})$$
 therefore, if $j$ (resp. $j'$) is an index appearing in the first (resp. second) direct sum decomposition, we have  $2j+2j'\geq 2n_{\ell}+1$. Whence the assertion.
\end{proof}

\subsubsection{A polarization}

For a nonnegative integer $\nu$, we define the primitive subspaces
\begin{equation}\nonumber
P_{r+\nu }:=\ker ((2i\pi N)^{\nu +1}:\gr^W_{r+\nu }H^{2*}_{\orb }(\mathcal{X})\rightarrow \gr^W_{r-\nu  -2}H^{2*}_{\orb}(\mathcal{X})).
\end{equation}
We now show that the Poincar\'e-Saito duality $g$ provides a polarization of the induced Hodge structure
\begin{equation}\label{eq:DecompPrim}
P_{r+\nu}=\oplus_p F^p P_{r+\nu} \cap \overline{F^{r+\nu -p} P_{r+\nu}}.
\end{equation}
First, and thanks to Lemma \ref{lemma:propertiesgorbifold}, it is possible to apply the construction recalled in the appendix: 
we define the nondegenerate bilinear form $g_{\nu}$ on $\gr^W_{r+\nu}H^{2*}_{\orb}(\mathcal{X})$ by putting 
\begin{equation}\nonumber
g_{\nu} (a , b) := g(\tilde{a} , (2i\pi N)^{\nu} \tilde{b})
\end{equation}
if $\tilde{a}, \tilde{b}\in W_{r+\nu}H^{2*}_{\orb}(\mathcal{X})$ represent $a,b\in \gr^W_{r+\nu}H^{2*}_{\orb}(\mathcal{X})$
and the sesquilinear form 
$$h_{\nu}  : P_{r+\nu}\times P_{r+\nu}\rightarrow \cit$$
by putting
\begin{equation}\label{eq:Sesquilinear}
h_{\nu} (u ,v) :=i^{2p-(r+\nu)} g_{\nu} (u , \overline{v})
\end{equation}
if $u\in F^p P_{r+\nu} \cap \overline{F^{r+\nu -p} P_{r+\nu}}$.

\begin{lemma}
\label{lemma:hHerm}
Let $\nu$ be a nonnegative integer.
\begin{enumerate}
\item The Hodge subspaces of the Hodge decomposition (\ref{eq:DecompPrim})
are orthogonal with respect to the sesquilinear form $h_{\nu}$.
\item The sesquilinear form $h_{\nu}$ is Hermitian.
\end{enumerate}
\end{lemma}
\begin{proof} 
By Lemma \ref{lemma:gF}, we have $g_{\nu} (F^p P_{r+\nu}, F^{r+\nu +1-p} P_{r+\nu})=0$ and, by Lemma \ref{lemma:propertiesgorbifold}, these orthogonality relations are preserved under conjugation. Therefore
$$h_{\nu} (F^p P_{r+\nu} \cap \overline{F^{r+\nu -p} P_{r+\nu}}, F^q P_{r+\nu} \cap \overline{F^{r+\nu -q} P_{r+\nu}})= 0$$
 if $r+\nu-q \geq r+\nu +1-p$ or if
$q\geq r+\nu+1-(r+\nu ) +p$ that is if $q\leq p-1$ or $q\geq p+1$. 
For the second point, we use the previous assertion, the reality of $g$ and the $(-1)^{r+\nu}$-symmetry of $g_{\nu}$.
\end{proof}

\begin{lemma} 
\label{lemma:Polar}
Let $\nu$ be a nonnegative integer. 
\begin{enumerate}
\item We have $P_{r+\nu}= \bigoplus_{\ell\in F} H^{n_{\ell}-\nu}_{\prim}(X_{\ell}):= \bigoplus_{\ell\in F}\ker (  N^{\nu+1}: H^{n_{\ell}-\nu}(X_{\ell})\rightarrow H^{n_{\ell}+\nu +2}(X_{\ell}))$.
\item Let $a\in H^{n_{\ell}-\nu}_{\prim}(X_{\ell})$. Then,
$a\in F^{r-[i_{\ell}]-k} P_{r+\nu}\cap \overline{F^{r+\nu-(r-[i_{\ell}]-k)} P_{r+\nu}}$ where $2k=n_{\ell}-\nu$.
\item  Let $a\in H^{n_{\ell}-\nu}_{\prim}(X_{\ell})$. Then,
\begin{equation} \label{eq:h}
h_{\nu} (a,a)=\frac{1}{(2\pi)^{r-\nu}}(-1)^{\alpha -i_{\ell}}\int_{\mathcal{X}_{\ell}}^{\orb} a \wedge N^{\nu} c(a)
\end{equation}
where $2(\alpha -i_{\ell})=2k=n_{\ell}-\nu$.
\end{enumerate}
In particular, $h_{\nu} (a,a)>0$ if $a\in F^p P_{r+\nu}\cap \overline{F^{r+\nu -p} P_{r+\nu}}$ and $a\neq 0$.
\end{lemma}
\begin{proof}
We keep the notations of the statement of the lemma.
\begin{enumerate}
\item By the very definition of the weight filtration (see Section \ref{sec:Weight}), we have
$$\gr^W_{r+\nu} H_{\orb}^{2*}= \bigoplus_{\ell\in F} H^{n_{\ell}-\nu} (X_{\ell})$$
whence the first assertion, because $N$ preserves each sector (see Remark \ref{rem:ComputeL}). 
\item By the very definition of the Hodge filtration (see Section \ref{sec:HodgeFiltration}), we have $a\in F^{r-[i_{\ell}]-k} P_{r+\nu}$. Using (\ref{eq:BasicConjugForm}), we get  
$\overline{a}\in  F^{r-[i_{\ell^{-1}}]-k} P_{r+\nu}$ therefore 
$$a\in \overline{F^{r-[i_{\ell^{-1}}]-k} P_{r+\nu}}= \overline{F^{r+\nu-(r-[i_{\ell}]-k)} P_{r+\nu}}$$
because $r-[i_{\ell^{-1}}]-k=\nu +[i_{\ell}]+k$ by (\ref{eq:FormuleDim}).
\item From the second point and the definition of $h_{\nu}$, we get (keeping in mind the definition of the Poincar\'e-Lefschetz conjugation (\ref{eq:BasicConjugForm})) 
\begin{equation}\nonumber 
h_{\nu} (a,a)=i^{[i_{\ell^{-1}}]-[i_{\ell}]}g(a, (2i\pi N)^{\nu} \overline{a})
=(2i\pi )^{\nu} i^{[i_{\ell^{-1}}]-[i_{\ell}]}g(a, N^{\nu} \overline{a})
\end{equation}
\begin{equation}\nonumber 
=(2i\pi )^{\nu -r} i^{[i_{\ell^{-1}}]-[i_{\ell}]} (-1)^{[\alpha]} (-1)^{\alpha-i_{\ell}}\int_{\mathcal{X}_{\ell}}^{\orb} a \wedge N^{\nu} c(a) = \frac{1}{(2\pi)^{r-\nu}}(-1)^{\alpha-i_{\ell}}\int_{\mathcal{X}_{\ell}}^{\orb} a \wedge N^{\nu} c(a)
\end{equation}
because $r-\nu =[i_{\ell}]+[i_{\ell^{-1}}]+2(\alpha-i_{\ell})$ thanks again to (\ref{eq:FormuleDim}).
\end{enumerate}
Now, the right hand term of (\ref{eq:h}) is positive if $a\neq 0$ (see for instance \cite[Theorem 2.2]{Fe}, \cite{Sai},   \cite{St}; the standard Hodge structure on $X_{\ell}$ is of Hodge-Tate type). 
The last assertion follows since, by the second point, $F^p P_{r+\nu}\cap \overline{F^{r+\nu -p} P_{r+\nu}}$ is a suitable direct sum of primitive vector spaces $H^{2k}_{\prim}(X_{\ell})$ and $h_{\nu}$ preserves this direct sum decomposition.
\end{proof}

\begin{theorem} \label{theo:Polar}
Assume that $H^{2*}_{\orb} (\mathcal{X})$ is equipped with the Poincar\'e-Lefschetz conjugation (\ref{eq:BasicConjugForm}) and let $g$ be the Poincar\'e-Saito orbifold duality.
Then, the tuple $$( H_{\orb }^{2*} (\mathcal{X}, \cit ),  H_{\orb ,\rit}^{2*}(\mathcal{X}),  F^{\bullet} ,  W_{\bullet}, N, g, r)$$
is a polarized mixed Hodge structure of weight $r$.
\end{theorem}
\begin{proof} 
Properties 1-5 of Definition \ref{def:PolarizedMHS} are now clear (Lemma \ref{lemma:propertiesgorbifold} provides the required properties for $g$). Property 6 follows from 
Lemma \ref{lemma:hHerm} and property 7 follows from Lemma \ref{lemma:Polar}.
\end{proof}

It is interesting to specialize the previous results to weighted projective spaces, for which we do not need to refer to any conceptual result (in particular the hard Lefschetz property) in order to get a polarization.

\begin{proposition} 
\label{prop:PolarWPS}
Let $P$ be a reduced simplex of weight $q=(q_0 ,\ldots , q_n )$ as in Example \ref{ex:Simplices} and let $\nu$ be a nonnegative integer. 
\begin{enumerate}
\item We have $P_{r+\nu}=\bigoplus_{\ell\in F,\ n_{\ell}=\nu} H^{0}(X_{\ell})$.
\item Let $a\in H^0 (X_{\ell})$ with $n_{\ell}=\nu$. Then $a\in F^{r-[i_{\ell}]} P_{r+\nu}\cap \overline{F^{r+\nu-(r-[i_{\ell}])} P_{r+\nu}}$.
\item Let $a\in H^0 (X_{\ell})$ with $n_{\ell}=\nu$. Then, 
\begin{equation} \label{eq:hWPS}
h_{\nu} (a,a)=\frac{1}{(2\pi)^{r-\nu}}\int_{\mathcal{X}_{\ell}}^{\orb} a \wedge N^{\nu} c(a).
\end{equation}
\end{enumerate}
In particular, the Poincar\'e-Saito duality provides a polarization of the mixed Hodge structure
$\MHS_{\orb}$.
\end{proposition}
\begin{proof}
Recall that
$\gr^W_{r+\nu} H_{\orb}^{2*}= \oplus_{\ell\in F} H^{n_{\ell}-\nu} (X_{\ell})$.
By \cite{Kawa} we have $\dim H^{2i} (X_{\ell}) =1$ and $H^{2i+1}(X_{\ell}) =0$ for all $\ell\in F$ therefore the map  
$$N^{\nu +1} : H^{n_{\ell}-\nu} (X_{\ell})\rightarrow H^{n_{\ell}+\nu +2} (X_{\ell})$$
 is injective if $\nu\leq n_{\ell}-2$ since $N^{n_{\ell}}\neq 0$ (see the remark below). This shows the first assertion. The remaining ones are shown as in Lemma \ref{lemma:Polar}.
\end{proof}

\begin{remark} 
\label{rem:CalculIntOrbiWPS}
In order to calculate the integrals (\ref{eq:hWPS}), we use \cite{Mann}: for instance, with the notations of loc. cit. we have $\omega:=\varphi^{-1} ([f])=\mu c_1^{\orb} (O_{\ppit (q)} (1))$ and
\begin{equation}\nonumber 
\int^{\orb}_{\ppit (q)} N^n =\mu^n \int_{\ppit (q)}^{\orb} (c_1^{\orb} (O_{\ppit (q)} (1)))^n =\frac{\mu^n}{q_0 q_1 \ldots q_n}
\end{equation}
by \cite[Proposition IV.3.13]{Mann} (this is the volume of the polar polytope $P^{*}$ of  $P$; recall that $\mu =q_0 +q_1 +\ldots +q_n$). Since the twisted sectors $\mathcal{X}_{\ell}$ are also weighted projective spaces, we have analogous formulas for $\int_{\mathcal{X}_{\ell}}^{\orb} N^{n_{\ell}}$: for $\ell \in F$ we have
\begin{equation}\nonumber 
\int^{\orb}_{\mathcal{X}_{\ell}} N^{n_{\ell}} =\frac{\mu^{n_{\ell}}}{\prod_{i\in I(\ell)} q_i}
\end{equation}
where $I(\ell):=\{i,\ q_i\ell\in\zit\}$. On the way, we get the hard Lefschetz isomorphisms
 \begin{equation}\nonumber 
N^{n_{\ell}-2k} : H^{2k}(X_{\ell})\stackrel{\sim}{\longrightarrow} H^{2(n_{\ell}-k)}(X_{\ell})
\end{equation}
for $2k\leq n_{\ell}$. 
\end{remark}

\begin{corollary} Let $P$ be a reduced simplex and let $\nu$ be a nonnegative integer. Then,
\begin{equation}\nonumber
P_{r+\nu}=\bigoplus_{p\in\zit} F^{p} P_{r+\nu}\cap \overline{ F^{r+\nu -p} P_{r+\nu}}
\end{equation} 
where 
\begin{equation}\nonumber
F^{p} P_{r+\nu}\cap \overline{ F^{r+\nu -p} P_{r+\nu}}= \bigoplus_{\ell\in F,\ n_{\ell}=\nu ,\ [i_{\ell}]=r-p} H^{0}(X_{\ell}).
\end{equation} 
\end{corollary}

\subsection{Alternate outtake: orbifold Poincar\'e-Fernandez duality}
\label{sec:Alternate}

In this section, we assume that $H^{2*}_{\orb} (\mathcal{X})$ is equipped with the Poincar\'e-Lefschetz conjugation (\ref{eq:BasicConjug}):
in this case,  and by Lemma \ref{lemma:Nreal}, $N$ is real while $2i\pi N$ is not, so we leave 
the singularity framework, but this setting 
may be of interest since it very close from the one considered by Fernandez in \cite{Fe} (see also Remark \ref{rem:KnownCases}) and goes along with the existence of K\"{a}hler packages on the orbifold cohomology, see Remark \ref{remark:Kahler Package}.
However, it should be emphasized that in this case we get a polarization only under the Hodge-Tate condition (\ref{eq:ConditionClef}).\\

For $0\leq \alpha\leq n$, we define the bilinear forms 
$$g_{\PF} : H^{2\alpha}_{\orb} (\mathcal{X})\times H^{2(n-\alpha )}_{\orb} (\mathcal{X})\rightarrow \cit$$
 as the direct sum over the sectors $X_{\ell}$ such that  $\alpha -i_{\ell}$ is a nonnegative integer of
\begin{equation}\nonumber 
g^{(\ell)}_{\PF} : H^{2(\alpha -i_{\ell})} (X_{\ell}) \times H^{2(n-\alpha-i_{\ell^{-1}})} (X_{\ell^{-1}})\rightarrow\cit 
\end{equation}
where now
\begin{equation} \label{def:OrbiBilinearPF}g^{(\ell)}_{\PF} (a,b)=(-1)^{\alpha -i_{\ell}}\int_{\mathcal{X}_{\ell}}^{\orb} a \wedge I^* (b).
\end{equation}
We will call it the {\em Poincar\'e-Fernandez duality}.
We have the counterpart of Lemma \ref{lemma:propertiesgorbifold} (its proof is analogous): the main difference is that $g_{\PF}$ is not $(-1)^r$-symmetric in general.

\begin{lemma} \label{lemma:propertiesgorbifoldBis}
The Poincar\'e-Fernandez duality satisfies the following properties:
\begin{enumerate}
\item $g_{\PF}$ is nondegenerate,
\item $g_{\PF}(b,a)=(-1)^{r-[i_{\ell^{-1}}]-[i_{\ell}]} g_{\PF}(a,b)$ if $a\in H^{2(\alpha -i_{\ell})} (X_{\ell})$ and $b\in H^{2(n-\alpha-i_{\ell^{-1}})} (X_{\ell^{-1}})$, 
\item $g_{\PF}$ is real with respect to the Poincar\'e-Lefschetz conjugation (\ref{eq:BasicConjug}),
\item $g_{\PF}(Na, b)+g_{\PF}(a, Nb)=0$.\qed 
\end{enumerate}
\end{lemma}

We define the sequilinear form 
$h_{\nu}^{\PF}$ on $P_{r+\nu}$ as in Section \ref{sec:OrbiPoincDual}, replacing $2i\pi N$ by $N$. Lemma \ref{lemma:gF} and Lemma \ref{lemma:hHerm} are still valid for $g_{\PF}$.

\begin{lemma}
\label{lemma:PositiveDefinitePrimitiveHTBis}
Assume that the Hodge-Tate condition (\ref{eq:ConditionClef}) holds. Then, the sesquilinear form 
$h_{\nu}^{\PF}$ is Hermitian and positive definite on the primitive space $P_{r+\nu}$.
\end{lemma}
\begin{proof}
Under the Hodge-Tate type assumption, and by Lemma \ref{lemma:propertiesgorbifoldBis},
the bilinear form $g_{\PF}$ is $(-1)^r$-symmetric therefore
$h_{\nu}^{\PF}$ is Hermitian as in Lemma \ref{lemma:hHerm} and the remaining positivity assertion is shown as in Theorem \ref{theo:Polar}.
\end{proof}

\begin{theorem} \label{theo:PolarAlternate}
Assume that $H^{2*}_{\orb} (\mathcal{X})$ is equipped with a Poincar\'e-Lefschetz conjugation (\ref{eq:BasicConjug}) and that the Hodge-Tate condition (\ref{eq:ConditionClef}) holds. 
Then,
the tuple $$( H_{\orb }^{2*} (\mathcal{X}, \cit ),  H_{\orb ,\rit}^{2*}(\mathcal{X}),  F^{\bullet} ,  W_{\bullet}, N, g_{\PF}, r)$$
is a polarized mixed Hodge structure of weight $r$.\qed 
\end{theorem}

\begin{proof}
The statement follows from Lemma \ref{lemma:propertiesgorbifoldBis} and Lemma \ref{lemma:PositiveDefinitePrimitiveHTBis}.
\end{proof}

\begin{remark}
\label{remark:Kahler Package}
Let us assume that the graded Milnor ring $A^*$ is integrally graded. 
Using the results of Section \ref{sec:HThpq} and Section \ref{sec:Alternate},
it is possible to show that it satisfies a K\"{a}hler package in the sense of 
\cite[Section 7]{AHK} (with respect to $[f]$, the class of $f$ in $A^1$) if and only if the Hodge-Tate condition (\ref{eq:ConditionClef}) holds (note that only the Hodge-Tate case is considered in loc. cit.). 
The corresponding Hodge-Riemann property should provide the log-concavity of some relevant sequences   
(see \cite{AHK}, \cite{Huh}, \cite{Stanley}, \cite{StaU}...). 
\end{remark}

 \appendix

\section{Appendix (polarized mixed-Hodge structures)}

\label{sec:PolarizedMHS}

 We recall the definition of polarized (mixed-)Hodge structures used in these notes.

\begin{definition}
\label{def:PolarizedHS}
A polarized Hodge structure of weight $r$ is a tuple $(H, H_{\rit},  F^{\bullet}, g, r)$ where
\begin{itemize}
\item $r$ is an integer,
\item $H_{\rit}$ is a finite dimensional $\rit$-vector space,
\item $H=H_{\rit}\otimes\cit$,
\item $g$ is a nondegenerate, $(-1)^r$-symmetric bilinear form on $H$, with real values on $H_{\rit}$, 
\item $F^{\bullet}$ is a decreasing and exhaustive filtration of $H$ (the Hodge filtration)
\end{itemize}
such that 
\begin{equation}\label{eq:PHS1}
	H=\oplus F^p \cap \overline{F^{r-p}},
\end{equation}
\begin{equation}\label{eq:PHS2}
g(F^p H, F^{r+1-p} H)=0,
\end{equation}
\begin{equation}\label{eq:PHS3}
i^{2p-r} g (u, \overline{u})>0\ \mbox{for all}\ u\in F^p H\cap \overline{F^{r -p} H},\ u\neq 0.
\end{equation}
We will also say that $g$ is a polarization of the pure Hodge structure
$(H, H_{\rit},  F^{\bullet})$ of weight $r$.
\end{definition}

\begin{remark} 
\label{rem:PolarizedPureHodge}
Let $h$ be the Hermitian form on $H$
defined by 
$$h (u, v):=i^{2p-r} g (u, \overline{v})$$
 if $u\in H^{p,r -p}:=F^p H\cap \overline{F^{r-p}H}$. By (\ref{eq:PHS2}) 
the Hodge subspaces $H^{p,r-p}$
are orthogonal with respect to $h$ (first Hodge-Riemann bilinear relations) and by (\ref{eq:PHS3}) $h$
is positive definite (second Hodge-Riemann bilinear relations).
\end{remark}

Let $N:V\rightarrow V$ be a linear map on a finite dimensional vector space over a field of characteristic $0$, which satisfies $N^{r+1}=0$.
Then (see \cite[Lemma 6.4]{Schm}),

\begin{itemize}
\item there exists a unique increasing filtration $W_{\bullet}$ of $V$  
\begin{equation}\nonumber
0\subset W_0 \subset W_1 \subset \cdots \subset W_{2r} =V
\end{equation}
such that 
$N(W_i )\subset W_{i-2}$ and 
\begin{equation}\label{eq:HLW}
N^{\ell} :  W_{r+\ell}/ W_{r+\ell -1} \stackrel{\sim}{\longrightarrow}  W_{r-\ell}/ W_{r-\ell -1} 
\end{equation}
This is the {\em weight filtration of $N$ centered at $r$}. 

\item We have the Lefschetz decomposition 
\begin{equation}\label{eq:LefDecGen}
\gr_{r+\ell}^W V=\oplus_{i\geq 0} N^i P_{r+\ell +2i}
\end{equation}
if $\ell \geq 0$ where the (primitive) subspaces $P_{r+\ell}$ of $\gr^W_{r+\ell}V$ are defined by
\begin{equation}\label{eq:PrimSubspaceDefGen}
P_{r+\ell}:=\ker (N^{\ell+1}:\gr^W_{r+\ell}V\rightarrow \gr^W_{r-\ell-2}V)
\end{equation}

if $\ell\geq 0$ and $P_{r+\ell}=0$ if $\ell <0$. 

\item If $N$ is an infinitesimal isometry of a $(-1)^r$-symmetric nondegenerate linear form $g$ on $V$, {\em i.e} 
\begin{equation}\label{eq:isom}
g(Na,b)+g(a,Nb)=0
\end{equation}
for all $a,b\in V$, then
$g(W_i , W_{i'})=0$ for $i+i' <2r$ and 
one can define a nondegenerate $(-1)^{r+\ell}$-symmetric bilinear form $g_{\ell}$ on $\gr^W_{r+\ell} V$ for $\ell\geq 0$ by
\begin{equation}\label{eq:gell}
g_{\ell} (a , b)= g(\tilde{a} , N^{\ell} \tilde{b})	
\end{equation}
where $\tilde{a}\in W_{r+\ell}$ ({\em resp.} $\tilde{b}\in W_{r+\ell}$) is a representative of $a$ ({\em resp.} $b$).
By (\ref{eq:isom}), the Lefschetz decomposition (\ref{eq:LefDecGen}) 
is orthogonal with respect to $g_{\ell}$. 
\end{itemize}
We apply the previous construction to the following situation: $V:=H_{\rit}$ is a finite dimensional $\rit$-vector space and $H$ is the complexification of $H_{\rit}$. 

\begin{definition}
\label{def:PolarizedMHS}
A polarized mixed Hodge structure of weight $r$
is a tuple $(H, H_{\rit},  F^{\bullet}, g, r, W_{\bullet}, N)$ such that:
\begin{enumerate}
\item $g$ is $(-1)^r$-symmetric nondegenerate bilinear form on $H$ with real values on $H_{\rit}$,
\item the nilpotent endomorphism $N$ of $H_{\rit}$ is an infinitesimal isometry of $g$,
\item $W_{\bullet}$ is the weight filtration of $N$ centered at $r$,
\item the filtration $F^{\bullet}\gr^W_m H$ defines a pure Hodge structure of weight $m$ on $\gr^W_m H$,
\item $N(F^p )\subset F^{p-1}$,
\item $g (F^p H, F^{r+1-p} H)=0$,
\item $i^{2p-(r+\ell )} g_{\ell} (u, \overline{u})>0$ for all $u\in F^p P_{r+\ell}\cap \overline{F^{r+\ell -p} P_{r+\ell}}$, $u\neq 0$.
\end{enumerate}
We will say that the tuple is pre-polarized if only the first six properties are satisfied.
\end{definition}

\begin{remark} 
\label{rem:PolarizedMixedHodge}
Let $(H, H_{\rit},  F^{\bullet}, g, r, W_{\bullet}, N)$ be a polarized mixed Hodge structure of weight $r$. 
The primitive subspace $P_{r+\ell}$ is a Hodge substructure as the kernel of a morphism of pure Hodge structures and
we get the (induced) Hodge decomposition
\begin{equation}\label{eq:DecompPrimGen}
P_{r+\ell} =\oplus_p F^p P_{r+\ell} \cap \overline{F^{r+\ell-p}P_{r+\ell}}.	
\end{equation}
This Hodge structure is polarized by $g_{\ell}$ in the sense of Remark \ref{rem:PolarizedPureHodge}.
\end{remark}

\section{Appendix: duality of $\mathcal{D}$-modules for weighted projective spaces}

\label{sec:DualityDmodulesGen}

In this appendix, we explain why the construction of the polarization in Section \ref{sec:OrbiPoincDual}
is inspired by the duality of $\mathcal{D}$-modules studied in \cite{Sab} and applied in \cite{DoSa2} to (mirror partners of) weighted projective spaces.

\subsection{Foreword: about the duality of $\mathcal{D}$-modules}

\label{sec:DualityDmodules}

 Given a $\cit [\tau ,\tau^{-1}]$-module $G$, we denote by $\overline{G}$ the $\cit$-vector space $G$ equipped with the new module structure $p(\tau).g:=p(-\tau )g$. 
We will denote by $\overline{g}$ the elements of $\overline{G}$.
If moreover $G$ is equipped with a connection $\partial_{\tau}$ then so is $G$ and we have
$\partial_{\tau} \overline{g}:= \overline{-\partial_{\tau} g}$. 
When $G$ is the Fourier-transform of the Gauss-Manin system of $f$ as in Section \ref{sec:SabbahMHS}, duality of $\mathcal{D}$-modules gives (see \cite[Section 3]{Sab} for details and also \cite[Section 2.7]{Sai1} for the local model) the existence of nondegenerate $\cit [\tau , \tau^{-1}]$-bilinear pairings 
\begin{equation}\nonumber
S: G \times \overline{G} \longrightarrow \cit [\tau ,\tau^{-1}]
\end{equation}
such that
\begin{enumerate}
\item $\tau \partial_{\tau} S( g' , \overline{g''})=S (\tau\partial_{\tau} g', \overline{g''})+
S (g', \overline{\tau\partial_{\tau} g''})$,
\item $S$ sends $G_0 \times \overline{G_0}$ in  $\cit [\theta ]\theta^{n}$,
\item $S (g'', \overline{g'})=(-1)^n \overline{S (g', \overline{g''})}$.	
\end{enumerate}
Such a pairing $S$
induces nondegenerate bilinear forms 
$$g_{\mathcal{D}}: H_{\alpha} \times H_{1-\alpha} \longrightarrow \cit \tau^{-1}\ \mbox{for}\ \alpha \neq 0\ \mbox{and}\ g_{\mathcal{D}}: H_{0} \times H_{0} \longrightarrow \cit.$$
This provides a $(-1)^{n-1}$-symmetric nondegenerate bilinear form $g_{\mathcal{D}}$ on $H_{\neq 0}$ and a $(-1)^n$-symmetric nondegenerate bilinear form on $H_0$ such that
$$g_{\mathcal{D}}(F^p , F^{r+1-p} )=0,$$
where $F^{\bullet}$ denotes the Hodge filtration of the mixed Hodge structure $\MHS_f$ (as always, $r=n$ on $H_0$ and $r=n-1$ otherwise).

\subsection{Application to weighted projective spaces}

\label{sec:DualityDmodulesWPS}

Let $P$ be a reduced simplex of weight $q=(q_0 ,\ldots ,q_n)$ in the sense of Example \ref{ex:Simplices} and let $f$ be the Laurent polynomial defined by 
$$f (u) =\sum_{b\in \mathcal{V}(P)} u^b$$
 as in (\ref{eq:fOrbi}). 
We first describe the pairing $g_{\mathcal{D}}$ on $H$. We follow closely \cite{DoSa2}.

Recall the classes $\omega_0 ,\ldots ,\omega_{\mu -1}$ in $G_0$ defined in \cite[Proposition 3.2]{DoSa2}, which give a Birkhoff normal form for the Brieskorn lattice $G_0$. We will denote by $\alpha_{i}$ the $V$-order of $\omega_{i}$. By \cite[Lemma 4.1]{DoSa2}, there exists a unique (up to a nonzero constant) pairing $S$ as above: it is given by the formulas (the coefficients of $f$ differs slightly from the one in {\em loc. cit.}, and this explains the mild modifications in the formulas defining $S$):

\begin{equation}\label{eq:PairingS}
S(\omega_k ,\overline{\omega_{\ell}})=\left\{ \begin{array}{ll}
c_{k,\ell} S(\omega_0 ,\overline{\omega_{n}}) & \mbox{if}\ 0\leq k\leq n\ \mbox{and}\ k+\ell=n\\
c_{k,\ell}\frac{1}{q^q}S(\omega_0 ,\overline{\omega_{n}}) & \mbox{if}\ n+1\leq k\leq \mu -1\ \mbox{and}\ k+\ell=\mu +n\\
0 & \mbox{otherwise}
\end{array}\right. 
\end{equation}

\noindent where $q^q =q_0 ^{q_0}\ldots q_n^{q_n}$ and $c_{k,\ell}$ is a positive rational number expressed with the help of the weights $q_0 ,\ldots , q_n$, see for instance \cite[Lemme VI.1.18]{Mann}.
Moreover, $S( \omega_k ,\overline{\omega_{\ell}})\in \cit \tau^{-n}$ and  
we will write 
$$S(\omega_0 ,\overline{\omega_{n}}):=S^o (\omega_0 , \omega_{n})\tau^{-n}.$$

The pairing $S$ induces a nondegenerate bilinear form  $g_{\mathcal{D}}$ on $H$ as follows (the brackets denote the class in $H$): 

\begin{itemize}
\item 
on $H_{\neq 0}$, 
\begin{equation}\label{eq:gH0}
g_{\mathcal{D}}( [\tau^{[\alpha_i]}\omega_i ] , [\tau^{[\alpha_j]}\omega_j ]) =(-1)^{[\alpha_j]} S(\omega_i , \overline{\omega_j})\tau^{n-1}
\end{equation}
 if $[\alpha_i ]+[\alpha_j ]=n-1$ and $0$ otherwise,
\item on $H_{0}$,  
\begin{equation}\label{eq:gneq0}
g_{\mathcal{D}}( [\tau^{\alpha_i}\omega_i ] , [\tau^{\alpha_j}\omega_j ]) =(-1)^{\alpha_j} S(\omega_i , \overline{\omega_j} )\tau^{n}
\end{equation}
 if $\alpha_i +\alpha_j =n$ and $0$ otherwise.
\end{itemize}
\noindent This bilinear form is $(-1)^{n-1}$-symmetric on $H_{\neq 0}$, $(-1)^{n}$-symmetric on $H_0$.

Let $k$ be such that $[\omega_{k}']:=[\tau^{[\alpha_{k}]}\omega_{k}]$ 
is a primitive element of $N$ (we will denote by $K$ the set of these indices): the space
$$B_{k}:=<N^j [\omega_{k}'], \ j=0,\ldots , \nu_{k}>$$
is a Jordan block of $N$ of size $\nu_{k}+1$ and we will call the integer $\nu_{k}$ is the {\em weight} of $[\omega_{k}']$. 
For $k \in K$ such that $k\geq n+1$, 
the class $[\omega_{\mu +n -k-\nu_k}']$ is also a primitive vector of $N$ of weight 
$\nu_{k}$:
we will put $\overline{k} :=\mu +n -k-\nu_k $ and $\overline{k+j}:= \overline{k} +j$ for $j=0,\ldots ,\nu_{k}$ (if $k\in [0,n]$ we simply define $\overline{k}:=k$).
Note that
\begin{equation}\label{eq:alphaalphabarre}
\alpha_{k}+\alpha_{\overline{k}}+\nu_{k}=n	
\end{equation}
because the $V$-order of $\omega_{\mu +n-k}$ is equal to $n-\alpha_{k}$ and also to 
$\alpha_{\overline{k}}+\nu_{k}$ (compare with (\ref{eq:FormuleDim})).
The next result is a motivation for the definition of Poincar\'e-Lefschetz conjugations (see Section \ref{sec:Conjug}).

\begin{lemma}\cite[Lemma 5.3]{DoSa2}
\label{lemma:ConjugWPSDoSa}
Let $k \in K$. Then,
\begin{equation}\label{eq:Conjug}
\overline{[\omega_{k}']}=\sum_{\ell =k}^{k +\nu_k}a_{k,\overline{\ell}} [\omega_{\overline{\ell}}']
\end{equation}
where $\overline{a_{k ,\overline{k}}} a_{\overline{k} ,k}=1$ (in particular, $a_{k ,\overline{k}}\neq 0$). \qed 
\end{lemma}

 In what follows, the double brackets denote the class in $\gr^W_{\bullet}H$ and
we consider the bilinear form $g_{\mathcal{D}, \ell}$ induced on $\gr^W_{r+\ell}$ by $g_{\mathcal{D}}$ as in (\ref{eq:gell}), using the nilpotent endomorphism $2i\pi N$. The two next results must be compared with Lemma \ref{lemma:Polar} and Theorem \ref{theo:Polar} respectively.

\begin{proposition}\label{prop:BoitePrimitive}
Let $[\omega_{k}']$ be a primitive element of $N$ of weight $\nu_{k_0}$. Then,
\begin{enumerate}
\item $[[\omega_{k}']]\in P_{r+\nu_{k_0}}$ and $P_{r+\nu_{k_0}}$ is generated by such classes,
\item 	$[[\omega_{k}']]\in F^{r-[\alpha_{k}]} P_{r+\nu_{k_0}}\cap \overline{F^{r+\nu_{k_0}-(r-[\alpha_{k}])} P_{r+\nu_{k_0}}}$,
\item for any primitive element $[\omega_{\ell}']$ of weight $\nu_{k_0}$,
\begin{equation}\label{eq:Formulag}
g_{\mathcal{D}, \nu_{k_0}} ([[\omega_{k}']],\overline{[[\omega_{\ell}']]})=\left\{ \begin{array}{ll}
0 & \mbox{if}\  \ell \neq k\\
a_{k ,\overline{k}} (2i\pi)^{\nu_{k_0}} (-1)^{r-[\alpha_{k}]} c_{k, \mu +n-k} \frac{1}{q^q} S^o (\omega_0 ,\omega_n) & \mbox{if}\ \ell =k \neq 0  \\
a_{0 , 0} (2i\pi)^{n} (-1)^{n} S^o (\omega_0 ,\omega_n) & \mbox{if}\ \ell=k =0
\end{array}\right.
\end{equation}
\item the Hodge subspaces $F^p P_{r+\nu_{k_0}} \cap \overline{F^{r+\nu_{k_0} -p} P_{r+\nu_{k_0}}}$ are orthogonal with respect to $g_{\mathcal{D}, \nu_{k_0}}( \bullet , \overline{\bullet})$.
\end{enumerate}
\end{proposition}
\begin{proof} 
The first point follows from the definition of the primitive subspaces.
 By the very definition
of the Hodge filtration, we have $[[\omega_{k}']]\in F^{r-[\alpha_{k}]} P_{r+\nu_{k_0}}$.
Now, we have also $[[\omega_{\overline{k}}']]\in P_{r+\nu_{k_0}}$ since 
$\nu_{k_0}= \nu_{\overline{k}_0}$ and
it follows from (\ref{eq:Conjug}) that 
$\overline{[[\omega_{k}']]}=a_{k,\overline{k}} [[\omega_{\overline{k}}']]\in  F^{r-[\alpha_{\overline{k}}]} P_{r+\nu_{k_0}}$. Therefore, 
$$[[\omega_{k}']]\in \overline{F^{r-[\alpha_{\overline{k}}]} P_{r+\nu_{k_0}}}= \overline{F^{r+\nu_{k_0}-(r-[\alpha_{k}])} P_{r+\nu_{k_0}}}$$
where the equality follows from (\ref{eq:alphaalphabarre}). Formulas (\ref{eq:Formulag}) follow from the equalities
\begin{equation}\nonumber 
g_{\mathcal{D}, \nu_{k_0}} ([[\omega_{k}']],\overline{[[\omega_{\ell}']]})
=a_{\ell ,\overline{\ell}} (2i\pi)^{\nu_{k_0}}g_{\mathcal{D}} ([\omega_{k}'], N^{\nu_{k_0}}[\omega_{\overline{\ell}}'])
=a_{\ell ,\overline{\ell}} (2i\pi)^{\nu_{k_0}}g_{\mathcal{D}} ([\omega_{k}'], [\omega_{\overline{\ell}+\nu_{k_0}}']),
\end{equation}
(\ref{eq:PairingS}), (\ref{eq:gH0}) and (\ref{eq:gneq0}).  The fourth point follows 
from 1, 2 and 3.
\end{proof}

Let $h_{\mathcal{D}, \ell}$ be the sesquilinear form defined on the primitive subspace $P_{r+\ell}$ by  
\begin{equation}\label{eq:hellBis}
 h_{\mathcal{D}, \ell} (u ,v) :=i^{2p-(r+\ell )} g_{\mathcal{D}, \ell} (u , \overline{v})
\end{equation}
if $u\in F^p P_{r+\ell} \cap \overline{F^{r+\ell -p} P_{r+\ell}}$.

\begin{proposition}
\label{prop:SesquilinearSimplex}
Let $[\omega_{k}']$ be a primitive element of $N$ of weight $\nu_{k_0}$. Then,
\begin{enumerate}
\item $h_{\mathcal{D}, \nu_{k_0}} ([[\omega_{k}']], [[\omega_{k}']])=\frac{1}{(2i\pi)^{r}}(2\pi)^{\nu_{k_0}+r} a_{k ,\overline{k}}c_{k, \mu +n-k}\frac{1}{q^q}S^o (\omega_0 ,\omega_n)$ if $k \neq 0$,
\item $h_{\mathcal{D}, n} ([[\omega_{0}']], [[\omega_{0}']])=\frac{1}{(2i\pi)^n}(2\pi)^{2n} a_{0 , 0}S^o (\omega_0 ,\omega_n)$,
\item $h_{\mathcal{D}, \nu_{k_0}} ([[\omega_{k}']],[[\omega_{\ell}']])=0$ for any primitive element $[\omega_{\ell}']$ of weight $\nu_{k_0}$ such that $\ell\neq k$.
\end{enumerate}
\end{proposition}
\begin{proof}
By Proposition \ref{prop:BoitePrimitive}, we have
$h_{\mathcal{D}, \nu_{k_0}} ([[\omega_{k}']],[[\omega_{k}']])=i^{2p-(r+\nu_{k_0})}  g_{\mathcal{D}, \nu_{k_0}} ([[\omega_{k}']],\overline{[[\omega_{k}']]})$
where $p=r-[\alpha_{k}]$
hence, using (\ref{eq:Formulag}),
\begin{equation}\nonumber
h_{\mathcal{D}, \nu_{k_0}} ([[\omega_{k}']], [[\omega_{k}']])=i^{2r-2[\alpha_{k}]-r} 
(-1)^{r-[\alpha_{k}]} (2\pi)^{\nu_{k_0}}
a_{k ,\overline{k}} c_{k, \mu +n-k}\frac{1}{q^q} S^o (\omega_0 , \omega_n )  
\end{equation}
\begin{equation}\nonumber 
=i^{-r}  (2\pi)^{\nu_{k_0}} a_{k ,\overline{k}} c_{k, \mu +n-k}\frac{1}{q^q} S^o (\omega_0 , \omega_n ).      
\end{equation}
This shows the first point, and the computations are analogous for the second one. 
The third point also follows from Proposition \ref{prop:BoitePrimitive}.
\end{proof}

\begin{remark}
Proposition \ref{prop:SesquilinearSimplex} shows that it is natural to consider the bilinear form $(2i\pi )^r g_{\mathcal{D}}$. 
The bilinear form $(2i\pi)^{-r} g_{\mathcal{D}}$ is basically the orbifold Poincar\'e-Saito duality defined in Section \ref{sec:OrbiPoincDual} (the Poincar\'e-Lefschetz conjugation (\ref{eq:BasicConjugForm}) corresponds to the case $a_{k ,\overline{k}}=1$ in Lemma 
\ref{lemma:ConjugWPSDoSa}): this shift from $(2i\pi )^r g_{\mathcal{D}}$ to $(2i\pi)^{-r} g_{\mathcal{D}}$ is due to the fact that we replace $(-1)^{r-[\alpha_k]}$ by $(-1)^{[\alpha_k]}$ in Formula (\ref{eq:Formulag}).
\end{remark}


\begin{thebibliography}{999} 
\bibitem[1]{AHK} Adiprasito, K., Huh, J., Katz, E.: {\em Hodge theory for combinatorial geometries}, Annals of Math., {\bf 188}, 2018, p. 381-452.
\bibitem[2]{BCS} Borisov, L., Chen, L.,  Smith, G.: {\em The orbifold Chow ring of toric Deligne-Mumford stacks}, 
J. Amer. Soc. , {\bf 18} (1), 2005, p. 193-215. 
\bibitem[3]{CCLT} Coates, T., Corti, A., Lee, Y., Tseng, H.: {\em The quantum orbifold cohomology of weighted projective spaces}, Acta Mathematica, {\bf 202}, 2009,  p. 139-193. 
\bibitem[4]{CKS} Cattani, E., Kaplan, A., Schmid, W.: {\em Degeneration of Hodge structures}, Annals of  Mathematics, {\bf 123} (3), 1986,  p. 457-535. 
\bibitem[5]{CLS} Cox, D., Little, J., Schenck, A.: {\em Toric varieties}, American Mathematical Society, {\bf 124}, 2010.
\bibitem[6]{Conrads} Conrads, H. : {\em Weighted projective spaces and reflexive simplices}, Manuscripta Math., {\bf 107}, 2002, p. 215-227.
\bibitem[7]{CR} Chen, W., Ruan, Y.: {\em A new cohomology theory for orbifolds},
Comm. Math. Phys., {\bf 248} (1), 2004, p. 1-31.
\bibitem[8]{Del} Deligne, P.: {\em Th\'eorie de Hodge : II}, Pub. Math. IHES, {\bf 40}, 1971, p. 5-57.
\bibitem[9]{D9} Douai, A.: {\em A canonical Frobenius structure}, Mat. Z., {\bf 261}, 2009, p. 625-648.
\bibitem[10]{D13} Douai, A.: {\em Hard Lefschetz properties and distribution of spectra in singularity theory and Ehrhart theory}, J. Singul., {\bf 23}, 2021, p. 116-126.                                                                                                                                                                                                                                                                                                                                                                                                                                                                                                                                                                                                                                                                                                                                                                                                                                                                                                                                                                                                                                                                                                                                                                                                                                                                                                                                                                                                                                                                                                                                                                                                                                                                                                                                                                                                                                                                                                                                                                                                                                                                                                                                                                                                                                                                                                                                                                                                                                                                                                                                                                                                                                                                                                                                                                                                                                                                                                                                                                                                                                                                                                                                                                                                                                                                                                                                                                                                                                                                                                                                                                                                                                                                                                                                                                                                                                                                                                                                                                                                                                                                                                                                                                      
\bibitem[11]{DConfin} Douai, A.: {\em From Hodge theory for tame functions to Ehrhart theory for polytopes}, Ann. Comb., {\bf 29} (4), 2025, p. 1077-1107. DOI:10.1007/s00026-024-00738-7
\bibitem[12]{DoMa} Douai, A., Mann, E.: {\em The small quantum cohomology of a weighted projective space, a mirror D-module and their classical limits}, Geom. Dedicata, {\bf 164}, 2013, p. 187-226.
\bibitem[13]{DoSa1} Douai, A., Sabbah, C.: {\em Gauss-Manin systems, Brieskorn
lattices and Frobenius structures I}, Ann. Inst. Fourier, {\bf 53} (4), 2003, p. 1055-1116.
 \bibitem[14]{DoSa2} Douai, A., Sabbah, C.: {\em Gauss-Manin systems, Brieskorn lattices and Frobenius structures II}, In : Frobenius Manifolds, C. Hertling and M. Marcolli (Eds.), Aspects of Mathematics E 36, 2004.
\bibitem[15]{Fe} Fernandez, J.: {\em Hodge structures for orbifold cohomology}, Proc. Amer. Math. Soc., {\bf 134} (9), 2006, p. 2511-2520.
\bibitem[16]{FK} Fleming, B., Karu, K.:  {\em Hard Lefschetz theorems for simple polytopes},  J. Alg. Comb., {\bf 32}, 2010, p. 227-239.       
\bibitem[17]{GKR} Gross, M., Katzarkov, L., Ruddat, H.:  {\em Towards mirror symmetry for varieties of general type},  Advances in Mathematic, {\bf 308}, 2017, p. 208-275.    
\bibitem[18]{HL} Harder, A., Lee, S.:  {\em Irregular Hodge numbers of stacky Clarke mirror pairs}, arXiv:2408.09016    
\bibitem[19]{Huh} Huh, J.:  {\em Combinatorial applications of the Hodge-Riemann relations},  Proc. Int. Cong. of Math 2018 Rio de Janeiro, {\bf 3}, p. 3079-3098.       
\bibitem[20]{Jia} Jiang, Y.:  {\em The Chen-Ruan cohomology of weighted projective spaces},  Canad. J. Math., {\bf 59}, 2007, p. 981-1007.
\bibitem[21]{KS} Katz, E., Stapledon, A.: {\em Tropical geometry, the motivic nearby  fiber and limit mixed Hodge numbers of hypersurfaces}. Res. Math. Sc. {\bf 3} (10), 2016. 
\bibitem[22]{Kawa} Kawasaki, T.: {\em Cohomology of twisted projective spaces and lens complexes} Math. Ann., {\bf 206}, 1973, p. 243-248.
\bibitem[23]{KKP} Katzarkov, L., Kontsevitch, M., Pantev, T.: {\em Bogomolov-Tian-Todorov theorems for Landau-Ginzburg models}, J. Differential Geometry, {\bf 105} (1), 2017, p. 55-117.
\bibitem[24]{K} Kouchnirenko, A.G.: {\em Poly\`edres de Newton et nombres de Milnor}, Invent. Math.,  {\bf 32}, 1976, p. 1-31.
\bibitem[25]{Mann} Mann, E.: {\em Cohomologie quantique orbifolde des espaces projectifs \`a poids}, arXiv:0510331.
\bibitem[26]{Mann1} Mann, E.: {\em Orbifold quantum cohomology of weighted projective spaces}, Journal of Algebraic Geometry, {\bf 17}, 2008, p. 137-166.
\bibitem[27]{MT} Matsui, Y., Takeuchi, K.: {\em Monodromy at infinity of polynomial maps and Newton polyhedra (with an appendix by C. Sabbah)}. Int. Math. Res. Not. IMRN, {\bf 8}, 2013, p. 1691-1746.
\bibitem[28]{McM} McMullen, P.: {\em On simple polytopes}, Inv. Math., {\bf 113}, 1993, p. 419-444.
\bibitem[29]{Payne} Payne, S.: {\em Ehrhart series and lattice triangulations}, Discrete Comput. Geom., {\bf 40} (3), 2008, p. 365-376.
\bibitem[30]{Ph} Pham, F.: {\em La descente des cols par les onglets de Lefschetz avec vue sur Gauss-Manin}. In ''Syst\`emes diff\'erentiels et singularit\'es'',  Ast\'erisque, Soc. Math. France,  {\bf 130}, 1985.
\bibitem[31]{Sab0} Sabbah, C.: {\em  Monodromy at infinity and Fourier transform}. Pub. RIMS, Kyoto Univ., {\bf 33}, 1997, p. 643-685. 
\bibitem[32]{Sab} Sabbah, C.: {\em Hypergeometric periods for a tame polynomial}, Portugalia Mathematicae, {\bf 63}, 2006, p. 173-226.
\bibitem[33]{Sab2} Sabbah, C.: {\em Some properties and applications of Brieskorn lattices}, Journal of Singularity theory, {\bf 18}, 2018, p. 238-247.
\bibitem[34]{SY} Sabbah, C., Yu, J.D: {\em On the irregular Hodge filtration of exponentially twisted mixed Hodge modules}, Forum Math. Sigma, {\bf 3}, 2015.
\bibitem[35]{SaiK} Saito, K.: {\em The higher residue pairings $K_F^{(k)}$ for a family of hypersurfaces singular points}, Singularities, Proc. of Symposia in Pure Math., American Mathematical Society, {\bf 40}, 1983, p. 441-463.
\bibitem[36]{Sai1} Saito, M.: {\em On the structure of Brieskorn lattices}, Ann. Inst. Fourier, {\bf 39}, 1989, p. 27-72.
\bibitem[37]{Sai} Saito, M.: {\em Mixed Hodge modules}, Pub. Res. Inst. Math. Sci., {\bf 26}, 1990, p. 221-333.
\bibitem[38]{ScSt} Scherk, J.,  Steenbrink, J. H. M.: {\em On the mixed Hodge structure on the cohomology of the Milnor fibre}, Math. Ann., {\bf 271}, 1985, p. 641-665.
\bibitem[39]{Schm} Schmid, W.: {\em Variation of Hodge structure: the singularities of the period mapping}, Invent. Math., {\bf 22}, 1973, p. 211-319.
\bibitem[40]{Stanley} Stanley, R. P.: {\em Combinatorial applications of the hard Lefschetz theorem}, in Proc. of the International Congress of Mathematicians (Warsaw, 1983), North-Holland, Amsterdam, 1984, p. 447-453.
\bibitem[41]{StaU} Stanley, R. P.: {\em Log-concave and unimodal sequences in algebra, combinatorics and geometry}, In Graph theory and its applications: East and West (Jinan, 1986), volume 576 of Ann. New York Acad. Sci, pages 500-535. New York Acad. Sci., New York, 1989.
\bibitem[42]{St} Steenbrink, J. H. M.: {\em Mixed Hodge structure on the vanishing cohomology}, Real and complex singularities (Proc. Ninth Nordic Summer School/NAVF Sympos. Math., Oslo 1976), Sijthoff and Noordhoff, Alphen aan den Rijn, 1977, p. 525-553.
\bibitem[43]{SZ} Steenbrink, J.H.M, Zucker, S.: {\em Variations of mixed Hodge structures I}, Inv. Math., {\bf 80}, 1985, p. 489-542. 
\bibitem[44]{Tak} Takeuchi, K.: {\em Geometric monodromies, mixed Hodge numbers of motivic Milnor fibers and Newton polyhedra}. arXiv:2308.09418
\bibitem[45]{Var} Varchenko, A. N.: {\em On the monodromy operator in vanishing cohomology and the operator of multiplication by $f$ in the local ring}, Sov. Math. Dokl., {\bf 24}, 1981, p. 248-252.
\bibitem[46]{Wa} Wang, H.: {\em Newton polyhedrons and Hodge numbers of non-degenerate Laurent polynomials}, arXiv:2403.02105.
\end{thebibliography}
\end{document}